\documentclass[reqno,12pt,a4]{amsart}
\usepackage{amssymb}
\usepackage{amsmath}

\usepackage[all]{xy}

\def\Hom{\mathop{\rm Hom}\nolimits}
\def\Ext{\mathop{\rm Ext}\nolimits}
\def\id{\mathop{\rm id}\nolimits}

\def\Ad{\mathop{\rm Ad}\nolimits}

\def\ev{\mathop{\rm ev}\nolimits}

\newtheorem{thm}{Theorem}[section]
\newtheorem{lem}[thm]{Lemma}

\newtheorem{rmk}[thm]{Remark}

\numberwithin{equation}{section}

\newcommand{\co}{\operatorname{co}}

\newcommand{\Bott}{\operatorname{Bott}}
\newcommand{\Aut}{\operatorname{Aut}}
\newcommand{\as}{\operatorname{as}}

\date{\today}
\author{V. Manuilov and K. Thomsen}
\title{Semi-invertible extensions and asymptotic homomorphisms}

\begin{document}

\maketitle

\begin{abstract}
We consider the semigroup $\Ext(A,B)$ of extensions of a separable
$C^*$-algebra $A$ by a stable $C^*$-algebra $B$ modulo unitary equivalence
and modulo asymptotically split extensions. This semigroup contains the
group $\Ext^{-1/2}(A,B)$ of invertible elements (i.e. of semi-invertible
extensions). We show that the functor $\Ext^{1/2}(A,B)$ is homotopy
invariant and that it coincides with the functor of homotopy classes of
asymptotic homomorphisms from $C(\mathbb T)\otimes A$ to $M(B)$ that map
$SA\subseteq C(\mathbb T)\otimes A$ into $B$.

\end{abstract}

\section{Introduction}

This is a study of a general structure in the extensions of a
separable $C^*$-algebra by another separable and stable $C^*$-algebra.
The significance of such extensions comes from many applications, but
is perhaps best illustrated by the fact that all the common homotopy
invariant and stable functors on the category of separable
$C^*$-algebras admit descriptions in terms of such $C^*$-extensions.
To explain our viewpoint on these extensions, which originates from
our work in \cite{MT3} and the problems which it naturally leads us
to consider, we must put the results and methods from \cite{MT3} into
perspective.

The main discovery in \cite{MT3} was that the E-theory of Connes
and Higson is the quotient of the unitary equivalence classes of
extensions by the asymptotically split extensions, provided the
$C^*$-algebras that play the roles of quotient and ideal in the
extensions are, respectively, suspended and stable. This reveals
that if the role of the split extensions, which has served as the
natural trivial extensions since the work of Brown, Douglas and
Fillmore, \cite{BDF1}, \cite{BDF2}, are replaced by the
asymptotically split extensions, then the question about
invertibility of the extensions disappear, at least when the
quotient is a suspended $C^*$-algebra. The significance of this
is stressed by the (albeit slowly) growing number of examples of
extensions which are not invertible in the BDF sense,
\cite{An}, \cite{Wass},\cite{S}, \cite{Ki}, \cite{Kapl}, \cite{HT} .
Among these, Kirchbergs examples are the most
striking in our optic because they show that the BDF semi-group
of extensions fail to be a group for a large class of naturally
occuring $C^*$-algebras, in cases where \emph{the homotopy classes}
of extensions {\it do} form a group.

An important point concerning the methods used in \cite{MT3} is that
they provide proofs of homotopy invariance in the class of unitary
equivalence classes of extensions modulo the asymptotically split
extensions by using the
relation to asymptotic homomorphisms given by the Connes-Higson
construction, \cite{CH}. This is a completely new approach to
homotopy invariance in the theory of $C^*$-algebra extensions
which is independent of the methods which were developed for
this in \cite{BDF2} and \cite{K}.

However, the methods in \cite{MT3} require in an essential way
that the $C^*$-algebra which plays the role of the quotient in
the extension is a suspended $C^*$-algebra. This is annoying
because it means that general $C^*$-algebra extensions must be
suspended in order to become amenable to the methods and results
of \cite{MT3}, and this is particularly frustrating because the
key tool from \cite{MT3}, the Connes-Higson construction, is
available for any $C^*$-extension. The most crucial reason for the
success of the methods developed in \cite{MT3} is that every extension
is semi-invertible, in the sense that it can be made asymptotically split
by adding another extension to it, when the quotient is a suspended $C^*$-algebra. One of the main questions left open by \cite{MT3} is therefore

{\bf Question:} Does the Connes-Higson construction, in the general
case, provide us with an isomorphism, from unitary equivalence
classes of semi-invertible extensions modulo asymptotically split extensions to
homotopy classes of asymptotic homomorphisms ?

At present we do not know, in the general case, if every extension
of a separable $C^*$-algebra by a separable stable $C^*$-algebra
is semi-invertible. All the examples
mentioned above of extensions that fail to be invertible in the BDF
sense may very well turn out to be semi-invertible. In fact, it
follows from \cite{MT3} that the examples of Kirchberg, \cite{Ki},
are semi-invertible. Thus we must also ask:

{\bf Question:} Are all extensions of a separable $C^*$-algebra
by a separable stable $C^*$-algebra semi-invertible ?

The main purpose here is to answer the first question by a qualified
'Yes'. More precisely we show that a variant of the Connes-Higson
construction, which takes the semi-invertibility of the extensions
into account, does give rise to an isomorphism. Unfortunately this
does not, in itself, answer the question for the genuine
Connes-Higson map.

\section{The group of semi-invertible extensions}

Let $A$ and $B$ be separable $C^*$-algebras, $B$ stable. Let
$M(B)$ be the multipler algebra of $B$ and $Q(B) = M(B)/B$ the
generalized Calkin-algebra of $B$. Let $q_{B} : M(B) \to Q(B)$
be the quotient map. The extensions of $A$ by $B$ will be
identified with $\Hom(A,Q(B))$; the $*$-homomorphisms from $A$ to
$Q(B)$. Two extensions $\varphi, \psi \in \Hom(A,Q(B))$ are
\emph{unitary equivalent} when there is a unitary $u \in M(B)$
such that $\Ad q_{B}(u)\circ \varphi = \psi$. An extension
$\psi \in \Hom(A,Q(B))$ is \emph{asymptotically split} when
there is an asymptotic homomorphism $\pi =(\pi_t)_{t \in
[1,\infty)} : A \to M(B)$ such that $q_{B} \circ \pi_t =
\psi$ for all $t$. Thanks to the stability of $B$, the
unitary equivalence classes in $\Hom(A,Q(B))$ form an abelian
semi-group: Choose isometries $V_1,V_2 \in M(B)$ such that
$V_1V_1^* + V_2V_2^* = 1$, and define $\varphi \oplus \psi$
to be the extension $a \mapsto \Ad q_{B} (V_1) \circ
\varphi(a) + \Ad q_{B} (V_2) \circ \psi(a)$. Then the
addition in the unitary equivalence classes in $\Hom(A,Q(B))$
is given by $[\varphi] + [\psi] = [\varphi \oplus \psi]$. This
addition is independent of the choice of isometries $V_1,V_2$,
subject to the condition that $V_1V_1^* + V_2V_2^* = 1$. We
say that the extension $\varphi \in \Hom(A,Q(B))$ is
\emph{semi-invertible} when there is another extension $\psi$
such that $\varphi \oplus \psi$ is asymptotically split. Both
the semi-invertible and the asymptotically split extensions
represent a semi-group in the unitary equivalence classes of
extensions; the latter contained in the first, and we denote
the 'quotient' by $\Ext^{-1/2}(A,B)$. Thus two semi-invertible
extensions, $\varphi$ and $\psi$, define the same element of
$\Ext^{-1/2}(A,B)$ if and only if there are asymptotically
split extensions, $\lambda_i, i =1,2$, such that $\varphi \oplus
\lambda_1$ is unitarily equivalent to $\psi \oplus \lambda_2$.
The main goal of the paper is to obtain a description of
$\Ext^{-1/2}(A,B)$ in terms of asymptotic homomorphisms. For this purpose we set $\Ext^{-1/2}(A,D) = \Ext^{-1/2}(A,D \otimes \mathbb K)$, where $\mathbb K$ is the 
$C^*$-algebra of compact operators on a separable infinite dimensional Hilbert space, when $D$ is a separable $C^*$-algebra which is not stable. Note that $\Ext^{-1/2}(A,D)$ is functorial (contravariantly) in an obvious way in the first variable 
$A$. In the second variable, $D$, there is a priori only functoriality with respect to quasi-unital $*$-homomorphisms, cf. \cite{H}, in a way we now describe. Given a quasi-unital $*$-homomorphism $\varphi : D \to D_1$, the tensor 
product, $\varphi \otimes \id_{\mathbb K} : D \otimes \mathbb K \to D_1 \otimes \mathbb K$, of $\varphi$ with the identity on $\mathbb K$ is again quasi-unital 
and admits therefore an extension $\widetilde{\varphi \otimes \id_{\mathbb K}}: M(D \otimes \mathbb K) \to M(D_1 \otimes \mathbb K)$ which, in turn, defines a 
$*$-homomorphism $\widehat{\varphi \otimes\id_{\mathbb K}} : Q(D \otimes \mathbb K) \to  Q(D_1 \otimes \mathbb K)$. 
We set $\varphi_*[\psi] = [\widehat{\varphi \otimes\id_{\mathbb K}} \circ \psi]$. When $e \in \mathbb K$ is a minimal non-zero projection, we define $s_D : D \to D \otimes \mathbb K$ by $s_D(d) = d \otimes e$.

\begin{lem}\label{gen1} ${s_D}_* : \Ext^{-1/2}(A,D) \to
\Ext^{-1/2}(A,D \otimes \mathbb K)$ is an isomorphism.
\end{lem}
\begin{proof} This is all very standard: As is well-known,
there is an isometry $V \in M(D \otimes \mathbb K \otimes \mathbb K)$
and an isomorphism $\gamma : D \otimes \mathbb K \to D \otimes
\mathbb K \otimes \mathbb K$ such that $\Ad V \circ \gamma =
s_D \otimes \id_{\mathbb K}$. It suffices therefore to show that
conjugation by the isometry $V$ induces the identity map of
$\Ext^{-1/2}(A, D \otimes \mathbb K)$, and this is clear because
conjugation by $V$ is just addition by the trivial
extension $0$.
\end{proof} \qed

In other words, the functor $\Ext^{-1/2}(A, -)$ is stable, and
there is no reason to distinquish between $\Ext^{-1/2}(A, B)$ and
$\Ext^{-1/2}(A, B \otimes \mathbb K)$ when $B$ is a stable separable
$C^*$-algebra.

\section{Pairing $\Ext^{-1/2}$ with $KK$-theory}\label{pairing}

In this section we prove homotopy invariance of $\Ext^{-1/2}$ in
the second variable. Homotopy invariance in the first variable is
an immediate consequence. Unlike the approach taken in \cite{MT3},
the proof hinges on Kasparov's homotopy invariance result from
\cite{K}, in the more abstract guise it was given by Higson in \cite{H}.

Recall that an asymptotic homomorphism $\varphi =
\left(\varphi_t\right)_{t \in [1, \infty)} : A \to B$ between
$C^*$-algebras is \emph{equi-continuous} when the family of maps,
$\varphi_t : A \to B, t \in [1,\infty)$, is an equi-continuous
family of maps. As is well-known any asymptotic homomorphism is
asymptotically equal to one which is equi-continuous. We shall make 
use of the following generalisation of this fact. The proof
exploits the so-called asymptotic algebra of a given $C^*$-algebra
$E$, via the Bartle-Graves selection theorem. Let $C_b\left([1,\infty),
E\right)$ be the $C^*$-algebra of continuous and norm-bounded
$E$-valued function on $[1,\infty)$ and $C_0\left([1,\infty), E\right)$
the ideal in $C_b\left([1,\infty), E\right)$ consisting of elements
$f$ for which $\lim_{t \to \infty} \|f(t)\| = 0$.
The \emph{asymptotic algebra} $\as(E)$ of $E$ is the quotient
$$
\as(E) = C_b\left([1,\infty), E\right)/C_0\left([1,\infty), E\right) .
$$

\begin{lem}\label{gen2} Let $A,B, D$ be $C^*$-algebras,
$A_0 \subseteq A$, $B_0 \subseteq B$ $C^*$-subalgebras and
$\chi : B \to D$ a $*$-homomorphism. Let
$\pi = \left(\pi_t\right)_{t \in [1,\infty)} :
A \to B$ be an asymptotic homomorphism and
$\mu : A \to D$ a $*$-homomorphism such that $\chi \circ \pi_t =
\mu$ for all $t \in [1,\infty)$. Assume that $\pi_t(A_0) \subseteq B_0$
for all $t \in [1,\infty)$.

It follows that there is an equi-continuous asymptotic homomorphism
$\hat{\pi} : A \to B$ such that
\begin{enumerate}
\item[1)] $\lim_{t \to \infty} \pi_t(a) - \hat{\pi}_t(a) = 0$ for all
$a \in A$,
\item[2)] $\chi \circ \hat{\pi}_t = \mu$ for all $t \in [1,\infty)$,
\item[3)] $\hat{\pi}_t(A_0) \subseteq B_0$ for all $t \in [1, \infty)$.
\end{enumerate}
\end{lem}
\begin{proof} Set $P = \{(a,b) \in A \oplus B: \ \mu(a) = \chi(b) \}$,
and $\tilde{\pi}_t(a) = \left(a, \pi_t(a)\right) \in P$. Then
$\tilde{\pi}$ is an asymptotic homomorphism, and defines in a
natural way a $*$-homomorphism $\overline{\pi} :
A \to \as(P)$
into the asymptotic algebra of $P$. Since $\pi_t(A_0) \subseteq B_0$
by assumption,
\begin{equation}\label{gene1}
\overline{\pi}\left(A_0\right) \subseteq  \as(P_0),
\end{equation}
where $P_0 = \{ (a,b) \in P: \ a \in A_0, \ b \in B_0 \}$.
It follows from the Bartle-Graves selection theorem that there is
a continuous lift $\Phi : A \to C_b\left( [1,\infty), P\right)$ of
$\overline{\pi}$. For a detailed account of the Bartle-Graves
selection theorem we refer to \cite{L}, where there is also an
important remark, Remark 2 on p. 114, that we shall use: Because
of (\ref{gene1}) we can choose $\Phi$ such that $\Phi(A_0) \subseteq P_0$. 
Set $\hat{\pi}_t(a) = p\left(\Phi(a)(t)\right)$, where
$p : P \to B$ is the projection to the second coordinate.
\end{proof} \qed

Let $A,B,D$ be separable $C^*$-algebras, $B$ and $D$ stable. Let
$\psi' \in \Hom(A,Q(B))$ be a semi-invertible extension, and $x$ an
element of $KK(B,D)$. $x$ is then represented, in the picture of
$KK$-theory obtained in \cite{T}, by a pair of $*$-homomorphisms
$\pi_{\pm} : M(B) \to M(D)$ such that
\begin{equation}\label{U3}
\pi_+(b)-\pi_-(b) \in D
\end{equation}
for all $b \in B$. Since $\psi'$ is semi-invertible, there is an
asymptotic homomorphism $\psi = (\psi_t)_{t \in [1,\infty)} :
A \to M_2(M(B))$, given in matrix notation as
$$
\psi_t = \left ( \begin{matrix} \psi^t_{11} & \psi^t_{12} \\
\psi^t_{21}  & \psi^t_{22} \end{matrix} \right),
$$
such that $q_{M_2(B)} \circ \psi_t$ is $t$-independent and
$q_B \circ \psi^t_{11} = \psi'$ for all $t \in [1,\infty)$.
In particular,
\begin{equation}\label{U4}
\psi^t_{12}(a), \psi^t_{21}(a)  \in B
\end{equation}
for all $t,a$. Note that by Lemma \ref{gen2} we can assume that
$\psi$ is equi-continuous. We will refer $\psi$ as a
\emph{trivialization} of $\psi'$. Set
$$
\left(\pi_{\pm} \times \psi\right)_t(a) =
q_{M_2(D)}\left ( \begin{matrix} \pi_+ \left(\psi^t_{11}(a)\right) &
\pi_+ \left(\psi^t_{12}(a)\right) \\ \pi_+ \left(\psi^t_{21}(a)\right)  &
\pi_- \left( \psi^t_{22}(a)\right) \end{matrix} \right) .
$$

\begin{lem}\label{U2} $ \left(\left(\pi_{\pm} \times
\psi\right)_t\right)_{t \in [1,\infty)} : A \to M_2(Q(D))$
is an asymptotic homomorphism.
\end{lem}
\begin{proof} Calculating modulo $M_2(D)$ we find that
\begin{equation*}
\begin{split}
&  \left(\begin{matrix} \pi_+ \left(\psi^t_{11}(a)\right) & \pi_+
\left(\psi^t_{12}(a)\right) \\ \pi_+ \left(\psi^t_{21}(a)\right)  &
\pi_- \left( \psi^t_{22}(a)\right) \end{matrix} \right)
\left ( \begin{matrix} \pi_+ \left(\psi^t_{11}(b)\right) & \pi_+
\left(\psi^t_{12}(b)\right) \\ \pi_+ \left(\psi^t_{21}(b)\right)  &
\pi_- \left( \psi^t_{22}(b)\right) \end{matrix} \right)    \\
& \\
& =\left(\begin{matrix} \pi_+ \left(\psi^t_{11}(a) \psi^t_{11}(b) +
\psi^t_{12}(a) \psi^t_{21}(b)\right) & \pi_+ \left(\psi^t_{11}(a)
\psi^t_{12}(b)\right) + \pi_+\left( \psi^t_{12}(a)\right) \pi_-
\left( \psi^t_{22}(b)\right) \\
\pi_+ \left(\psi^t_{21}(a)\psi^t_{11}(b) \right) +
\pi_-\left(\psi^t_{22}(a)\right)\pi_+\left( \psi^t_{21}(b)\right)  &
\pi_+ \left( \psi^t_{21}(a) \psi^t_{12}(b)\right) +\pi_-
\left(\psi^t_{22}(a)\psi^t_{22}(b)\right) \end{matrix} \right) \\
&\\
& =    \left(\begin{matrix} \pi_+ \left(\psi^t_{11}(a)
\psi^t_{11}(b) + \psi^t_{12}(a) \psi^t_{21}(b)\right) &
\pi_+ \left(\psi^t_{11}(a) \psi^t_{12}(b) +  \psi^t_{12}(a)
\psi^t_{22}(b)\right) \\ \pi_+ \left(\psi^t_{21}(a)\psi^t_{11}(b) +
\psi^t_{22}(a) \psi^t_{21}(b)\right)  & \pi_- \left( \psi^t_{21}(a)
\psi^t_{12}(b) + \psi^t_{22}(a)\psi^t_{22}(b)\right) \end{matrix} \right),
 \end{split}
\end{equation*}
thanks to (\ref{U3}) and (\ref{U4}). Since $\psi$ is an asymptotic
homomorphism the last expression is asymptotically equal to
\begin{equation*}
 \left(\begin{matrix} \pi_+ \left(\psi^t_{11}(ab)\right) & \pi_+
\left(\psi^t_{12}(ab)\right) \\ \pi_+ \left(\psi^t_{21}(ab)\right)  &
\pi_- \left( \psi^t_{22}(ab)\right) \end{matrix} \right).
\end{equation*}
\end{proof} \qed

Let $\hat{\Theta} : M_2(Q(D)) \to Q(D)$ be a $*$-isomorphism induced
by two isometries $V_1,V_2 \in M(D)$ such that $V_1V_1^* + V_2V_2^* = 1$.

\begin{lem}\label{U6}  Assume that $\psi' \in \Hom(A,Q(B))$ is asymptotically 
split. It follows that there are
asymptotic homomorphisms $\mu, \nu : A \to M(D)$ such that
$$
 \left(\hat{\Theta} \circ \left(\pi_{\pm} \times \psi\right)_t\right)
\oplus \left(q_D \circ \mu_t\right) = \left( q_D \circ \nu_t\right)
$$
for all $t$.
\end{lem}
\begin{proof} Since $\psi'$ is asymptotically split, 
there is an asymptotic homomorphism $\theta : A \to M(B)$
such that $\psi^t_{11}(a) - \theta_t(a) \in B$ for all $t,a$.
By Lemma \ref{gen2} we can assume that $\theta$ is equi-continuous. Let
$$
E_1 = \pi_-(B),
$$
$$ \mathcal F = \left\{ \pi_- \circ \theta_t(a) + \mathbb C \pi_-(1) :
\ t \in [1,\infty), a \in A\right\},
$$
and
$$
X = \left\{ \theta_t(a) + \psi^t_{11}(a') + \mathbb C 1:
\ t \in [1,\infty), a,a' \in A \right\}.
$$
Both $X$ and $\mathcal F$ are separable sets since the involved
asymptotic homomorphisms are equi-continuous. Let $E_2$ be the
$C^*$-algebra generated by $\left\{ \pi_-(x) - \pi_+(x) :
\ x \in X \right\}$. Then $E_1E_2 \subseteq D$ and
$\mathcal F E_1 \subseteq E_1$. Thus Kasparov's technical theorem,
\cite{K}, provides us with $M,N \in M(D)$ such that
$M,N \geq 0, M+ N = 1, ME_1 \subseteq D, N E_2 \subseteq D$ and
$[N,\mathcal F] \subseteq D$. Set
$$
U = q_{M_2(D)} \left( \begin{matrix} \pi_-(1) \sqrt{M}  &
- \pi_-(1)\sqrt{N} \\ \pi_-(1)\sqrt{N} & \pi_+(1) \sqrt{M} \end{matrix}
\right) .
$$
Since $\pi_-(1) \in \mathcal F$, we see that both $N$ and $M$
commute with $\pi_-(1)$ modulo $D$. Since $1 \in X$ and
$N E_2 \subseteq D$, we see that $N\pi_-(1) = N \pi_+(1)$ modulo $D$.
In particular, $N$ and $M$ both commute with $\pi_-(1)$ and
$\pi_+(1)$ modulo $D$. It follows that $U + U^* \geq 0$.  Using
this, and that $UU^* = U^*U = q_{M_2(D)}\left( \begin{smallmatrix}
\pi_-(1) & 0 \\ 0 & \pi_+(1) \end{smallmatrix} \right)$, we can
lift $U$ to an element $V \in M_2(M(D))$ such that $V + V^* \geq 0$ and
$VV^* = V^*V =  \left( \begin{smallmatrix} \pi_-(1) & 0 \\ 0 & \pi_+(1)
\end{smallmatrix} \right)$. Set $S =  \left( \begin{smallmatrix} 0 & 1 \\
1 & 0 \end{smallmatrix} \right) V$, and note that
$$
S = \left( \begin{matrix} \pi_-(1) \sqrt{N}  &  \pi_+(1)\sqrt{M} \\
\pi_-(1)\sqrt{M} & -\pi_-(1) \sqrt{N} \end{matrix} \right)
\quad \mbox{modulo}\ M_2(D).
$$
It follows that
$$
W = \left( \begin{smallmatrix} S & {} \\ {}  &
1 \end{smallmatrix} \right) +
\left( \begin{smallmatrix}  0 & 1 - \pi_+(1) & {}   \\
1 - \pi_-(1) & 0 & {} \\ {}  & {} & 0 \end{smallmatrix} \right)
$$
is a unitary in $M_3(M(D))$ such that
$$
W =  \left( \begin{smallmatrix}  \pi_-(1) \sqrt{N}  & \pi_+(1)\sqrt{M}  &
0   \\ \pi_-(1)\sqrt{M} &  - \pi_-(1) \sqrt{N} & 0 \\ 0  & 0 & 1
\end{smallmatrix} \right)   + \left( \begin{smallmatrix}  0 &
1 - \pi_+(1) & {}   \\ 1 - \pi_-(1) & 0 & {} \\ {}  & {} & 0
\end{smallmatrix} \right)
\quad\mbox{modulo}\ M_3(D).
$$
It follows from the properties of $N$ and $M$ that
\begin{equation*}
W  \left( \begin{matrix}  \pi_-\left(\theta_t(a)\right)   & 0  & 0   \\
0  &  \pi_+\left( \psi^t_{11}(a)\right) &   \pi_+
\left( \psi^t_{12}(a)\right) \\ 0  & \pi_+\left(\psi^t_{21}(a)\right) &
\pi_-\left( \psi^t_{22}(a)\right) \end{matrix} \right) =
\left( \begin{matrix}  \pi_+\left(\theta_t(a)\right)   & 0  & 0   \\
0  &  \pi_-\left( \psi^t_{11}(a)\right) &
\pi_-\left( \psi^t_{12}(a)\right) \\ 0  &
\pi_-\left(\psi^t_{21}(a)\right) &  \pi_-\left( \psi^t_{22}(a)\right)
\end{matrix} \right)W
\end{equation*}
modulo $M_3(D)$ for all $t,a$.
Set $\mu_t = q_D \circ \pi_- \circ \theta_t$ and
$$
\nu^0_t =    \left( \begin{matrix}  \pi_+\left(\theta_t(a)\right)   &
0  & 0   \\  0  &  \pi_-\left( \psi^t_{11}(a)\right) &
\pi_-\left( \psi^t_{12}(a)\right) \\ 0  &
\pi_-\left(\psi^t_{21}(a)\right) &  \pi_-\left( \psi^t_{22}(a)\right)
\end{matrix} \right) .
$$
Finally, we choose an appropriate isomorphism $\Theta_0 : M_3(M(D))
\to M(D)$ and set $\nu_t = \Theta_0 \circ \nu^0_t$.
\end{proof} \qed

\begin{lem}\label{U7} Let $\varphi', \psi' \in \Hom(A,Q(D))$ be
semi-invertible extensions with trivializations $\varphi$ and
$\psi$, respectively. There is then a trivialization $\lambda$ of
$\varphi' \oplus \psi'$ such that $\pi_{\pm} \times \lambda$ is
unitarily equivalent to  $\left(\pi_{\pm} \times \psi \right)
\oplus \left(\pi_{\pm} \times \varphi\right)$.
\end{lem}
\begin{proof}
Let $V_1,V_2 \in M(B)$ be the isometries used to define the addition
in $\Ext^{-1/2}(A,B)$. Then $\psi' \oplus \varphi' = q_D \circ
\left( \Ad V_1 \circ \psi^t_{11} + \Ad V_2 \circ \varphi^t_{11}\right)$
for all $t$, and
$$
\lambda_t =  \left( \begin{matrix} \Ad V_1 \circ \psi^t_{11} +
\Ad V_2 \circ \varphi^t_{11}  &  \Ad V_1 \circ \psi^t_{12} +
\Ad V_2 \circ \varphi^t_{12} \\  \Ad V_1 \circ \psi^t_{21} +
\Ad V_2 \circ \varphi^t_{21}  &  \Ad V_1 \circ \psi^t_{22} +
\Ad V_2 \circ \varphi^t_{22} \end{matrix} \right)
$$
is a trivialization of $\psi' \oplus \varphi'$. Note that
\begin{equation*}
\begin{split}
&\left( \pi_{\pm} \times \lambda \right)_t = \\
& \left( \begin{matrix} \Ad \pi_+(V_1) \circ \pi_+ \circ \psi^t_{11} +
\Ad \pi_+(V_2) \circ \pi_+ \circ \varphi^t_{11}  &  \Ad \pi_+(V_1)
\circ \pi_+ \circ \psi^t_{12} + \Ad \pi_+(V_2)  \circ \pi_+ \circ
\varphi^t_{12} \\  \Ad \pi_+(V_1) \circ \pi_+ \circ \psi^t_{21} +
\Ad \pi_+(V_2) \circ \pi_+ \circ \varphi^t_{21}  &  \Ad \pi_-(V_1)
\circ \pi_- \circ \psi^t_{22} + \Ad \pi_-(V_2) \circ \pi_- \circ
\varphi^t_{22} \end{matrix} \right)  .
\end{split}
\end{equation*}
Modulo $D$ we have that
\begin{equation*}
\begin{split}
& \pi_+(V_1) \left(\pi_+ \circ \psi^t_{12}(a)\right) \pi_-(V_1^*)  =
\pi_+\left(V_1\psi^t_{12}(a)\right) \pi_-\left(V_1^*\right) \\
& =  \pi_-\left(V_1\psi^t_{12}(a)\right) \pi_-\left(V_1^*\right)
\ \ \ \ \ \ \text{(by (\ref{U3}) and (\ref{U4}))} \\
& =\pi_-\left(V_1\psi^t_{12}(a)V_1^*\right) \\
& =\pi_+\left(V_1\psi^t_{12}(a)V_1^*\right)  \ \ \ \ \ \
\text{(by (\ref{U3}) and (\ref{U4}))} \\
& = \Ad \pi_+(V_1) \circ \pi_+ \circ \psi^t_{12}(a).
\end{split}
\end{equation*}
Via similar considerations regarding $\Ad \pi_+(V_2)  \circ \pi_+
\circ \varphi^t_{12},  \Ad \pi_+(V_1) \circ \pi_+ \circ \psi^t_{21}$,
and $\Ad \pi_+(V_2) \circ \pi_+ \circ \varphi^t_{21}$, we see that
\begin{equation*}
\begin{split}
&\left( \pi_{\pm} \times \lambda \right)_t = \\
& \Ad \left( \begin{matrix} \pi_+(V_1)    & 0  \\ 0  & \pi_-(V_1)
\end{matrix} \right) \circ \left( \begin{matrix}
\pi_+ \circ \psi^t_{11}   & \pi_+ \circ \psi^t_{12}  \\
\pi_+ \circ \psi^t_{21}  & \pi_- \circ \psi^t_{22}
\end{matrix} \right) +  \Ad \left( \begin{matrix} \pi_+(V_2)    & 0  \\
0  & \pi_-(V_2)  \end{matrix} \right) \circ \left( \begin{matrix}
\pi_+ \circ \varphi^t_{11}   & \pi_+ \circ \varphi^t_{12}  \\
\pi_+ \circ \varphi^t_{21}  & \pi_- \circ \varphi^t_{22}
\end{matrix} \right)  ,
\end{split}
\end{equation*}
modulo $M_2(D)$. We conclude that
\begin{equation*}
\left(\pi_{\pm} \times \lambda \right)_t = \Ad S_1 \circ
\left( \pi_{\pm} \times \psi \right)_t + \Ad S_2 \circ
\left( \pi_{\pm} \times \varphi\right)_t
\end{equation*}
modulo $M_2(D)$, where
$$
S_i = \left( \begin{smallmatrix} \pi_+(V_i) & {} \\ {} & \pi_-(V_i)
\end{smallmatrix} \right) + \left( \begin{smallmatrix}
1- \pi_+(1) & {} \\ {} & 1 -\pi_-(1) \end{smallmatrix} \right) ,
$$
$i = 1,2$, are isometries in $M_2(M(D))$ such that
$S_1S_1^* + S_2S_2^* = 1$. Thus, up to unitary equivalence, we have
that $\pi_{\pm} \times \lambda  = \left(\pi_{\pm} \times \psi \right)
\oplus \left(\pi_{\pm} \times \varphi \right)$.
\end{proof} \qed

We now introduce \emph{the basic construction} of \cite{MT5}.
Given an equi-continuous asymptotic homomorphism $\varphi  =
(\varphi_t)_{t \in [1,\infty)} : A \to Q(D)$, the basic
construction gives us a genuine extension $\varphi^f \in
\Hom (A,Q(D))$. The construction goes as follows: Let $b$ be a strictly
positive element in $D$ of norm $\leq 1$. A {\it unit sequence} (cf.
\cite{MT3})
in $D$ is a sequence $\{u_n\}_{n=0}^{\infty} \subseteq
E$ such that
\begin{enumerate}
\item[u1)] there is a continuous function $f_n : [0,1] \to [0,1]$ which
is zero in a neighbourhood of $0$ and $u_n = f_n(b)$,
\item[u2)] $u_{n+1}u_n = u_n $ for all $n$,
\item[u3)] $\lim_{n \to \infty} u_n x = x$ for any $x \in D$.
\end{enumerate}
Unit sequences exist by elementary spectral theory. Given a unit
sequence $\{u_n\}$ we set $\Delta_0 = \sqrt{u_0}$ and $\Delta_j =
\sqrt{u_j - u_{j-1}}, j \geq 1$. Note that u2) implies that
\begin{equation}\label{BB}
\Delta_i\Delta_j = 0, \ |i-j|\geq 2.
\end{equation}
Let $\hat{\varphi}_t : A \to M(D)$ be an equi-continuous lift of
$\varphi$, cf. Lemma 2.1 of \cite{MT5}. There exists a sequence
$t_1 < t_2 < t_3 < \dots$ in $[1,\infty)$ such that
$\{\varphi_{t_n}\}_{n = 1}^{\infty}$ is a discretization of $\varphi$ and
\begin{enumerate}
\item[t1)] $\lim_{n \to \infty} \sup_{t \in [t_n,t_{n+1}]}
\|\hat{\varphi}_t(a) - \hat{\varphi}_{t_n}(a)\| = 0$ for all $a \in A$, and
\item[t2)] $t_n \leq n$ for all $n \in \mathbb N$,
\end{enumerate}
cf. Lemma 3.3 of \cite{MT5}. We say that the pair $\left( \left(\hat{\varphi}_t\right)_{t \in [1,\infty)},
\{u_n\}_{n = 0}^{\infty}\right)$
is a \emph{compatible pair} for $\varphi$ when $
\lim_{n \to \infty} \sup_{t \in [1, n+1]}  \|u_n\hat{\varphi}_t(a) -
\hat{\varphi}_{t}(a)u_n\|  = 0$ for 
all $a \in D$. Compatible pairs exist by Lemma 3.2 of \cite{MT5}.
Given such a pair, there is a sequence $n_0 < n_1 < n_2 < \dots $ in
$\mathbb N$ such that
$$
n_i - n_{i-1} > i + 1
$$
for all $i \geq 1$,
$$\lim_{i \to \infty} \sup_{j \geq n_i} \sup_{t \in [1,i+3]}
\left( \|\left(1-u_j\right)\left( \hat{\varphi}_t(a)\hat{\varphi}_t(b) -
\hat{\varphi}_t(ab)\right)\| -  \|\varphi_t(a)\varphi_t(b) -
\varphi_t(ab)\| \right) = 0,
$$
\begin{equation*}
\lim_{i \to \infty} \sup_{j \geq n_i} \sup_{t \in [1, i+3]}
\left( \|\left(1 -u_j\right) \left( \hat{\varphi}_t(a) +
\lambda \hat{\varphi}_t(b) - \hat{\varphi}_t(a + \lambda b)\right)\| -
\|\varphi_t(a) + \lambda \varphi_t(b) - \varphi_t(a + \lambda b)\|
\right) = 0,
\end{equation*}
$$ \lim_{i \to \infty} \sup_{j \geq n_i} \sup_{t \in [1, i+3]}
\left( \|\left(1-u_j\right)\left(\hat{\varphi}_t(a^*) -
\hat{\varphi}_t(a)^*\right)\| -  \|\varphi_t(a^*) -
\varphi_t(a)^*\|\right) = 0
$$
for all $a,b \in A$ and all $\lambda \in \mathbb C$, cf.
Lemma 3.4 of \cite{MT5}. The quadruple
$$\left( \left(\hat{\varphi}_t\right)_{t \in [1,\infty)},
\{u_n\}_{n = 0}^{\infty}, \{n_i\}_{i=0}^{\infty},
\{t_i\}_{i=1}^{\infty}\right)
$$
is called \emph{the folding data}. Given the folding data there
is then an extension $\varphi^f : A \to Q(D)$ such that
$$
\varphi^f(a) = q_D\left(\sum_{j=0}^{\infty} \Delta_j
\hat{\varphi}_{t_{j+1}}(a) \Delta_j\right)
$$
for all $a \in A$, cf. Lemma 3.5 of \cite{MT5}. We will
refer to $\varphi^f$ as \emph{a folding} of $\varphi$.

We claim that we can define a map
$$
\Ext^{-1/2}(A,B)  \ni [\psi'] \mapsto \pi_{\pm}
\bullet [\psi'] \in \Ext^{-1/2}(A,D)
$$
by setting $\pi_{\pm} \bullet [\psi']  =
\left[\left( \hat{\Theta} \circ \left(\pi_{\pm}
\times \psi\right) \right)^f\right]$, where
$\psi$ is an arbitrary trivialization of $\psi'$.
To see that this recipe is well-defined we must show
that $ \left[\left(\hat{\Theta} \circ \left(\pi_{\pm}
\times \psi\right) \right)^f\right]$ is independent of
all the choices involved in its construction, and in
fact only depends on the class $[\psi']$ of $\psi' \in
\Ext^{-1/2}(A,B)$. For the first purpose, let $\varphi'
\in \Hom(A,Q(B))$ be an extension such that $\varphi'
\oplus \psi'$ is asymptotically split, and let $\varphi$
and $\psi$ be trivializations of $\varphi'$ and $\psi'$,
respectively. It follows then from Lemma \ref{U6} and
Lemma \ref{U7} that there are asymptotic homomorphisms
$\nu^1,\nu^2 : A \to M(D)$ such that
$$
\left(\hat{\Theta} \circ \left(\pi_{\pm} \times \varphi
\right)_t \right) \oplus \left( \hat{\Theta} \circ
\left(\pi_{\pm} \times \psi \right)_t \right) \oplus q_D
\circ \nu^1_t = q_D \circ \nu^2_t
$$
for all $t \in [1,\infty)$. Then Lemma 4.4 of \cite{MT5} implies that
$$
\left[\left( \hat{\Theta} \circ \left(\pi_{\pm} \times
\psi\right) \right)^f\right] = - \left[\left(\left( \hat{\Theta}
\circ \left( \pi_{\pm} \times \varphi \right)\right) \oplus
\left(q_D \circ \nu^1\right)\right)^f\right]  = -
\left[ \left( \hat{\Theta} \circ \left( \pi_{\pm} \times
\varphi \right) \right)^f\right]
$$
in $\Ext^{-1/2}(A,D)$. Thus $\left[ \left( \hat{\Theta}
\circ \left(\pi_{\pm} \times \psi\right) \right)^f\right]$
only depends on $\psi' \in \Hom(A,Q(D))$. To show that
$ \left[ \left(\hat{\Theta} \circ \left(\pi_{\pm} \times \psi'\right)
\right)^f\right]$ only depends on the class of $\psi'$ in
$\Ext^{-1/2}(A,B)$, it suffices now, thanks to Lemma \ref{U7}
and Lemma 4.5 of \cite{MT5}, only to show that the class is
not changed when $\psi'$ is replaced by an extension unitarily
equivalent to it. We leave this to the reader.

We want to show that the map $\pi_{\pm} \bullet - : \Ext^{-1/2}(A,B)
\to \Ext^{-1/2}(A,D)$ only depends on the class of $\left(\pi_+,
\pi_-\right)$ in $KK(B,D)$. For this purpose, we shall use the
homotopy-invariance theorem of Higson, cf. Section III of \cite{H}.

Thanks to Lemma \ref{gen1}, our construction above gives rise to a pairing
\begin{equation}\label{U11}
\Ext^{-1/2}\left( A, C \otimes B  \right) \to \Ext^{-1/2}\left( A, C\right)
\end{equation}
with quasi-unital Fredholm modules for $B$ in the sense of \cite{H}
for all separable $C^*$-algebras $A,C$ and $B$. This goes as follows:
Let $\varphi_{\pm} : B \to M(\mathbb K)$ be a quasi-unital Fredholm
pair in the sense of \cite{H}, i.e. $\varphi_{\pm}$ are quasi-unital
$*$-homomorphisms such that $\varphi_+(b) - \varphi_-(b) \in \mathbb K$
for all $b \in B$. Then $\id_C \otimes \varphi_{\pm} :
C \otimes B \to C \otimes M(\mathbb K) \subseteq M(C \otimes
\mathbb K)$ are quasi-unital $*$-homomorphisms, and admit canonical
extensions $\overline{\id_C \otimes \varphi_{\pm}} : M(C \otimes B)
\to M(C \otimes \mathbb K)$. Note that $\overline{\id_C \otimes
\varphi_{+}}(y) - \overline{\id_C \otimes \varphi_{-}}(y) \in
C \otimes \mathbb K$ for all $y \in C \otimes B$. Hence
$ x \mapsto {{s_C}_*}^{-1} \left(\overline{\id_C \otimes \varphi_{\pm}}
\bullet x \right)$ is a homomorphism , defining the desired
pairing (\ref{U11}) with quasi-unital Fredholm modules. By the quasi-unital
version of Higson's result, Theorem 3.1.4 and the remarks in the
first paragraph of Section 3.3 in \cite{H}, homotopy invariance
of $\Ext^{-1/2}(A,-)$ will now follow if we can show that the
pairing constructed above has the following properties
(cf. 3.1.3a -- 3.1.3f of \cite{H}):
\begin{enumerate}
\item[a)] $ h^*\left(\pi_{\pm} \bullet [\psi']\right) =
\left(\tilde{h} \circ \pi_{\pm}\right) \bullet [\psi']$,
when $h : D \to D'$ is a quasi-unital $*$-homomorphism with
canonical extension $\tilde{h} : M(D) \to M(D')$;
\item[b)] $\pi_{\pm} \bullet [\psi'] + \varphi_{\pm} \bullet
[\psi'] = (\pi_+,\varphi_-) \bullet [\psi']$, when $\pi_- = \varphi_+$;
\item[c)] $ ( \pi_{\pm} \oplus \pi) \bullet [\psi'] =  \pi_{\pm}
\bullet [\psi']$ for every $*$-homomorphism $\pi : M(B) \to M(D)$;
\item[d)] $\pi_{\pm} \bullet [\psi']=[\psi']$ when $\pi_+=\id:M({\mathbb
K})\to M({\mathbb K})$ and $\pi_-=0$;
\item[e)] $ \left(\Ad U \circ \pi_{\pm}\right) \bullet [\psi'] =
\pi_{\pm} \bullet [\psi']$ when $U \in M(D)$ is a unitary;
\item[f)] $( \pi, \Ad V \circ \pi) \times [\psi'] = 0$, when
$\pi : M(B) \to M(D)$ is a $*$-homomorphism and $V \in M(D)$ is
a unitary such that $V = 1$ modulo $D$.
\end{enumerate}
Of these, d) and f) are trivial and c) and e) follow from Lemma 4.4 and
Lemma 4.5 of \cite{MT5}, respectively. To prove b), we use Kasparov's
technical theorem in the following way: Let $E_2 \subseteq M_2(M(D))$
be the $C^*$-algebra generated by elements of the form
$$
\left( \begin{matrix} \pi_+ \left( \psi^t_{11}(a)\right) & \pi_+
\left( \psi^t_{12}(a)\right) \\ \pi_+\left( \psi^t_{21}(a) \right) &
0 \end{matrix} \right) \ \text{or} \ \left( \begin{matrix} \pi_-
\left( \psi^t_{11}(a)\right) & \pi_-\left( \psi^t_{12}(a)\right) \\
\pi_-\left( \psi^t_{21}(a) \right) & 0 \end{matrix} \right),
$$
$t \in [1,\infty), a \in A$, and $\mathcal F \subseteq M_2(M(D))$
the subspace spanned by elements of the form
$$
\left( \begin{matrix} \pi_+\left( \psi^t_{11}(a)\right) &  \pi_+
\left( \psi^t_{12}(a)\right) \\  \pi_+\left( \psi^t_{21}(a)\right) &
\pi_-\left( \psi^t_{22}(a)\right) \end{matrix} \right) \ \text{or}
\  \left( \begin{matrix} \varphi_+\left( \psi^t_{11}(a)\right) &
\varphi_+\left( \psi^t_{12}(a)\right) \\  \varphi_+
\left( \psi^t_{21}(a)\right) &  \varphi_-\left( \psi^t_{22}(a)\right)
\end{matrix} \right),
$$
$t \in [1,\infty), a \in A$. Let $E \subseteq M(D)$ be the
$C^*$-subalgebra consisting of the elements $m \in M(D)$ with
the property that $\pi_+(b)m , m \pi_+(b) \in D$ for all $b \in B$.
Since $\pi_-(b) = \pi_+(b) = \varphi_+(b) = \varphi_-(b)$ modulo
$D$ when $b \in B$, we might as well have used $\pi_-, \varphi_-$
or $\varphi_+$ instead of $\pi_+$ to define $E$. Note that
$ \pi_-\left( \psi^t_{22}(a)\right) - \varphi_-
\left( \psi^t_{22}(a)\right) \in E$ for all $a \in A$ and all
$t$, and that
$$
\mathcal F \left( \begin{smallmatrix} D & D \\ D & E \end{smallmatrix}
\right) \cup   \left( \begin{smallmatrix} D & D \\ D & E
\end{smallmatrix} \right) \mathcal F \subseteq
\left( \begin{smallmatrix} D & D \\ D & E \end{smallmatrix} \right) .
$$
We can therefore choose a separable $C^*$-subalgebra $E_1$ of
$\left( \begin{smallmatrix} D & D \\ D & E \end{smallmatrix} \right)$
containing
$$
\left( \begin{matrix} 0 & 0 \\ 0 &  \pi_-\left( \psi^t_{22}(a)\right) -
\varphi_-\left( \psi^t_{22}(a)\right) \end{matrix} \right)
$$
for all $a \in A$ and all $t$ such that $\left[\mathcal F, E_1\right]
\subseteq E_1$. Note that $E_1 E_2 \subseteq M_2(D)$. Kasparov's
technical theorem gives us elements $0 \leq N, M \in M_2(M(D))$
such that $N+ M  =1$, $[M,\mathcal F] \subseteq M_2(D)$, $M E_1
\subseteq M_2(D)$ and $N E_2 \subseteq M_2(D)$. Then
$$
U = \left( \begin{smallmatrix} - \sqrt{M} & \sqrt{N} \\ \sqrt{N} &
\sqrt{M} \end{smallmatrix} \right)
$$
is a unitary in $M_4(M(D))$ with the property that
\begin{equation*}
\begin{split}
&U \left( \begin{matrix} \pi_+\left( \psi^t_{11}(a)\right) &
\pi_+\left( \psi^t_{12}(a)\right) & 0 & 0 \\  \pi_+
\left( \psi^t_{21}(a)\right) & \pi_-\left( \psi^t_{22}(a)\right) &
0 & 0 \\  0 & 0 & \varphi_+\left( \psi^t_{11}(a)\right) &
\varphi_+\left( \psi^t_{12}(a)\right)  \\
0 & 0 & \varphi_+\left( \psi^t_{21}(a)\right) &
\varphi_-\left( \psi^t_{22}(a)\right) \end{matrix} \right) U^* \\
& =  \left( \begin{matrix} \pi_+\left( \psi^t_{11}(a)\right) &
\pi_+\left( \psi^t_{12}(a)\right) & 0 & 0 \\
\pi_+\left( \psi^t_{21}(a)\right) & \varphi_-\left( \psi^t_{22}(a)\right) &
0 & 0 \\  0 & 0 & \varphi_+\left( \psi^t_{11}(a)\right) &
\varphi_+\left( \psi^t_{12}(a)\right)  \\ 0 & 0 &
\varphi_+\left( \psi^t_{21}(a)\right) &
\varphi_+\left( \psi^t_{22}(a)\right)  \end{matrix} \right) ,
\end{split}
\end{equation*}
modulo $M_4(D)$, for all $t$ and $a$. This shows that
$\left(\pi_{\pm} \times \psi\right) \oplus \left( \varphi_{\pm}
\times \psi\right)$ is unitarily equivalent to
$\left( \left( \pi_+,\varphi_-\right) \times \psi\right) \oplus
\left(q_D \circ \nu\right)$, where $\nu : A \to M(D)$ is an
asymptotic homomorphism. It follows then from Lemma 4.4 and
Lemma 4.5 of \cite{MT5} that $\left[ \left(\hat{\Theta} \circ
\left(\pi_{\pm} \times \psi\right) \right)^f\right] +
\left[\left( \hat{\Theta} \circ \left( \varphi_{\pm}
\times \psi\right) \right)^f\right] =
\left[ \left( \hat{\Theta} \circ
\left( \left( \pi_+,\varphi_-\right) \times
\psi\right) \right)^f\right]$ in $\Ext^{-1/2}(A,D)$, proving b).
To prove a), note first that
$\hat{\Theta'} \circ \left( \left(\tilde{h} \circ \pi_{\pm}\right)
\times \psi \right)_t =  \Ad q_{D'}(U) \circ \hat{h} \circ \hat{\Theta}
\circ \left( \pi_{\pm} \times \psi\right)_t$ for all $t$, where
$U \in M(D')$ is the unitary $U =W_1 \left[\tilde{h}(V_1^*) +
1 - \tilde{h}(1)\right] + W_2 \tilde{h}(V_2^*)$, when $V_1,V_2$
and $W_1,W_2$ are the isometries used to define $\Theta$ and $\Theta'$,
respectively, and $\hat{h}:Q(D)\to Q(D)$ is induced by $\tilde{h}$.
Thanks to Lemma 4.5 of \cite{MT5} it remains therefore
only to prove that $\left[\left(  \hat{h} \circ \left( \hat{\Theta}
\circ \left( \pi_{\pm} \times \psi\right)\right)\right)^f\right] =
\left[\hat{h} \circ \left(  \hat{\Theta} \circ \left( \pi_{\pm}
\times \psi\right)  \right)^f\right]$, or if we set $\varphi =
\hat{\Theta} \circ \left( \pi_{\pm} \times \psi\right)$, that
$\left[ \hat{h} \circ \varphi^f\right] = \left[\left( \hat{h}
\circ \varphi\right)^f\right]$ in $\Ext^{-1/2}(A,D')$. Let
$\chi : A \to Q(D)$ be an equi-continuous asymptotic homomorphism
with the property that $\varphi \oplus \chi$ asymptotically splits.
Since $\tilde{h}$ is strictly continuous on norm-bounded sets,
we see that
$$
\hat{h} \circ \varphi^f(a) = q_{D'}\left( \sum_{j=0}^{\infty}
\tilde{h}\left(\Delta_j\right)
\tilde{h}\left( \hat{\varphi}_{t_{j+1}}(a)\right)\tilde{h}
\left( \Delta_j\right)\right) ,
$$
for a given tuple of folding data
$\left( \left(\hat{\varphi}_t\right)_{t \in [1,\infty)},
\{u_n\}_{n = 0}^{\infty}, \{n_i\}_{i=0}^{\infty},
\{t_i\}_{i=1}^{\infty}\right)$. Then the proof of
Lemma 4.4 in \cite{MT5} shows that
$$
\left[\hat{h} \circ \varphi^f\right] = -
\left[ \hat{h} \circ \chi^f\right]
$$
in $\Ext^{-1/2}(A,D')$. Since $ \left[  \hat{h} \circ \chi^f\right] =
-  \left[ \left( \tilde{h} \circ \varphi\right)^f \right]$ by
Lemma 4.4 of \cite{MT5}, we obtain the desired conclusion.

For any $C^*$-algebra $E$ we denote in the following the
$C^*$-algebra $C[0,1] \otimes E$ by $IE$. It follows that the
functor $\Ext^{-1/2}(A, -)$ is homotopy invariant, in the sense
that the point evaluations $\pi_t : IB  \to B, t \in [0,1]$,
induce the same maps ${\pi_t}_* : \Ext^{-1/2}(A, IB) \to
\Ext^{-1/2}(A,B)$ for any separable $C^*$-algebra $B$. When
this is established it is easy to make $\Ext^{-1/2}(A, - )$
functorial with respect to arbitrary $*$-homomorphisms; if
$h : B \to B_1$ is a $*$-homomorphism it follows from \cite{T}
that $h \otimes \id_{\mathbb K} : B \otimes \mathbb K \to
B_1 \otimes \mathbb K$ is homotopic to a quasi-unital
$*$-homomorphism $g : B \otimes \mathbb K \to B_1 \otimes
\mathbb K$, unique up to homopy, and we set $h_*[\psi] =
[\hat{g} \circ \psi]$, when $\psi \in \Hom(A,Q(B \otimes \mathbb K))$.
Thus we have obtained the following.

\begin{thm}\label{Thm1} For every separable $C^*$-algebra $A$,
$\Ext^{-1/2}(A, - )$ is a homotopy invariant functor, from the
category of separable $C^*$-algebras to the category of abelian groups.
\end{thm}

It follows also from the homotopy invariance that the pairing
$\pi_{\pm} \bullet - : \Ext^{-1/2}(A,B) \to \Ext^{-1/2}(A,D)$
constructed above only depends on the class of $\pi_{\pm}$ in
$KK(B,D)$. Thus we have in fact a pairing
$$
KK(B,C) \times \Ext^{-1/2}(A, B) \to  \Ext^{-1/2}(A, C)
$$
for all separable $C^*$-algebras $A,B$ and $C$.

It follows from Theorem \ref{Thm1} that two semi-invertible
extensions of $A$ by $B \otimes \mathbb K$ define the same
element of $\Ext^{-1/2}(A,B)$ if and only if they are homotopic
via a semi-invertible homotopy. Specifically, two semi-invertible
extensions $\varphi, \psi : A \to Q(B \otimes \mathbb K)$, define
the same element of $\Ext^{-1/2}(A,B)$ if and only if there is a
semi-invertible extension $\Phi : A \to Q(IB \otimes \mathbb K)$
such that $\hat{\pi_0} \circ \Phi = \psi$ and $\hat{\pi_1} \circ
\Phi = \varphi$, where $\hat{\pi_i}: Q(IB\otimes \mathbb K) \to
Q(B\otimes \mathbb K), i = 0,1$, are the $*$-homomorphisms
induced by the point evaluations $\pi_0,\pi_1 : IB \otimes
\mathbb K \to B \otimes \mathbb K$. It is this consequence
of Theorem \ref{Thm1} that we shall make intensive use of
in the following. But let us point out that the homotopy
invariance of $\Ext^{-1/2}(A,B)$ in the second variable,
$B$, implies the homotopy invariance in the first variable.

\begin{thm}\label{Thm2} For every separable $C^*$-algebra $B$,
$\Ext^{-1/2}(-, B)$ is a homotopy invariant functor, from the
category of separable $C^*$-algebras to the category of abelian groups.
\end{thm}
\begin{proof} Let $\varphi, \psi  : A \to D$ be homotopic
$*$-homomorphisms between separable $C^*$-algebras. Thus there
is a $*$-homomorphism $\Phi : A \to ID$ such that $\pi_0 \circ
\Phi = \varphi, \pi_1 \circ \Phi = \psi$. Let $\chi : D \to
Q(B \otimes \mathbb K)$ be a semi-invertible extension. Let
$\tau : IQ(B \otimes \mathbb K) \to Q(IB\otimes \mathbb K)$
be the canonical inclusion. Then
$$
\hat{\pi_0} \circ   \tau \circ \left(\id_I \otimes \chi\right)
\circ \Phi = \varphi
$$
and
$$
\hat{\pi_1} \circ   \tau \circ \left(\id_I \otimes \chi\right)
\circ \Phi = \psi,
$$
so $[\psi] = {\pi_1}_*\left[ \tau \circ \left(\id_I \otimes
\chi\right) \circ \Phi\right] =  {\pi_0}_*\left[ \tau \circ
\left(\id_I \otimes \chi\right) \circ \Phi\right] = [\varphi]$
in $\Ext^{-1/2}(A,B)$ by Theorem \ref{Thm2}.
\end{proof} \qed

\section{Extended asymptotic homomorphisms}

In this section $A$ and $B$ are separable $C^*$-algebras. Let
$J \subseteq A$ be a $C^*$-subalgebra of $A$. An asymptotic
homomorphism $\varphi = \left(\varphi_t\right)_{t \in [1,\infty)} :
A \to M(B)$ is \emph{extended from $J$} when $\varphi_t(J)
\subseteq B$ for all $t \in [1,\infty)$. If the context identifies
the subalgebra $J$, we say simply that $\varphi$ is \emph{extended}.
If $\varphi$ is extended from $J$ and $q_B \circ \varphi_t =
q_B \circ \varphi_1$ for all $t$, or equivalently,
$\varphi_t(x)-\varphi_1(x) \in B$ for all $x \in  A$ and all $t$,
we say that $\varphi$ is \emph{constantly extended from $J$} or
just \emph{constantly extended}. Two (constantly) extended
asymptotic homomorphisms $\varphi, \psi : A \to M(B)$ are
\emph{homotopic} when there is an (constantly) extended asymptotic
homomorphism $\Phi : A \to M(IB)$ such that $\widetilde{\pi_0}
\circ \Phi_t = \varphi_t$ and $\widetilde{\pi_1} \circ \Phi_t =
\psi_t$ for all $t \in [1,\infty)$, where $\widetilde
{\pi_s} : M(IB ) \to M(B)$ is the $*$-homomorphism induced by the
point evaluation $\pi_s : IB \to B$. Homotopy is an equivalence
relation in both cases, and we denote by $[[A,J;B]]$ the homotopy
classes of extended asymptotic homomorphisms, and by $[\{A,J;B\}]$
the homotopy classes of constantly extended asymptotic homomorphisms.
The set $[[A,J;B]]$ has been introduced and studied in \cite{Guentner} in relation
to relative E-theory.

\begin{thm}\label{Thm3} The canonical (forgetful) map $[\{A,J;B\}]
\to [[A,J;B]]$ is a bijection.

\end{thm}

\begin{proof} Surjectivity: Let $\varphi :A \to M(B)$ be an
extended asymptotic homomorphism. To show that $\varphi$ is
homotopic to a constantly extended asymptotic homomorphism we
may assume that $\varphi$ is equi-continuous since it is
asymptotically identical, and hence homotopic, to such an
extended asymptotic homomorphism by Lemma \ref{gen2}. By
Lemma 4.1 of \cite{MT5} there is a continuous increasing
function $r : [1,\infty) \to [1,\infty)$ such that
$\left(\varphi_{r(t)}\right)_{t \in [1,\infty)}$ is
uniformly continuous, in the sense that the function
$t \mapsto \varphi_{r(t)}(a)$ is uniformly continuous
for all $a \in A$. Since $\psi =
\left(\varphi_{r(t)}\right)_{t \in [1,\infty)}$ is
homotopic to $\varphi$, it suffices to show that $\psi$
is homotopic to a constantly extended asymptotic homomorphism.
Let $F_1 \subseteq F_2 \subseteq F_3 \subseteq \dots $ be a
sequence of finite sets with dense union in $A$, such that
$\bigcup_n F_n \cap J$ is dense in $J$. Let $\epsilon_0 \geq
\epsilon_1 \geq \epsilon_2 \geq   \dots$ be a sequence in
$]0,1[$ chosen so small that
\begin{equation}\label{1u}
\left\|\left[a, \psi_t(x)\right]\right\| \leq \epsilon_i,
\ \left\|\left[b, \psi_t(x)\right]\right\| \leq \epsilon_{i}
\Rightarrow \left\|\left[\sqrt{b-a},\psi_t(x)\right] \right\|
\leq 2^{-i -1},
\end{equation}
for all $t \in [1,i+2], \ x \in F_i$, and
\begin{equation}\label{2'}
\left\| a\psi_t(x) -         \psi_t(x) \right\| \leq \epsilon_i,
\left\| b\psi_t(x) -   \psi_t(x) \right\| \leq \epsilon_{i}
\Rightarrow \left\|\sqrt{b-a} \psi_t(x)\right\| \leq 2^{-i-1},
\end{equation}
for all $t \in [1,i+2], \ x \in F_i \cap J$, when $0 \leq a
\leq b \leq 1$. Let $v_0,v_1,v_2, \dots$, be a unit sequence in
$B$ such that
\begin{equation}\label{1}
\left\|\left[v_i, \psi_t(x)\right]\right\| \leq \epsilon_i,
\ t \in [1,i+2], \ x \in F_i,
\end{equation}
\begin{equation}\label{2}
\left\| v_i\psi_t(x) -       \psi_t(x) \right\| \leq \epsilon_i,
\ t \in [1,i+2], \ x \in F_i \cap J.
\end{equation}

Let $n_0 < n_1 < n_2 < \dots$ be a sequence in $\mathbb N$ such
that $n_i - n_{i-1} > i +1 $ for all $i \geq 1$. We claim that
there are continuous paths $u_i(t), t \in [1,\infty),
i = 0,1,2, \dots $, in $B$, such that $u_0(t) \leq u_1(t)
\leq u_2(t) \leq \dots$ is a unit sequence in $B$ for all $t$,
\begin{equation}\label{3'}
u_i(1) = v_{n_i},
\end{equation}
and for $t \in [n,n+1]$ one has
\begin{equation}\label{3}
u_i(t) \in \co \{v_j : j \geq n\}
\end{equation}
for all $i = 0,1,2, \dots$,
and
\begin{equation}\label{4}
u_i(t) = u_{i}(1),\ i \geq n+1.
\end{equation}
(In particular, at integer points we have
the following equations:
 $$
\begin{array}{llllll}
u_0(1)=v_{n_0}, & u_0(2)=v_{n_1}, & u_0(3)=v_{n_2}, & u_0(4)=v_{n_3},
& u_0(5)=v_{n_4}, & \ldots;\\
u_1(1)=v_{n_1}, & u_1(2)=v_{n_1+1}, & u_1(3)=v_{n_2+1}, & u_1(4)=v_{n_3+1},
& u_1(5)=v_{n_4+1}, & \ldots;\\
u_2(1)=v_{n_2}, & u_2(2)=v_{n_2}, & u_2(3)=v_{n_2+2}, & u_2(4)=v_{n_3+2},
& u_2(5)=v_{n_4+2}, & \ldots;\\
u_3(1)=v_{n_3}, & u_3(2)=v_{n_3}, & u_3(3)=v_{n_3}, & u_3(4)=v_{n_3+3},
& u_3(5)=v_{n_4+3}, & \ldots;
\end{array}
 $$
etc.) The construction is the same as the
construction of $\{w_i(t)\}_{i=0}^{\infty}$ in the proof of
Lemma 4.4 in \cite{MT5}: Assume that $\{u_i(t)\}_{i=0}^{\infty}$,
$t \in [1,k]$, have been constructed, and that
$v_{n_{k-1}} \leq u_0(k) \leq u_k(k) = v_{n_k}$.
We construct then $\{u_i(t)\}_{i=0}^{\infty}$, $t \in [k,k+1]$,
as follows. Since $n_{k+1} - n_k > k+1$, we have that
$$
v_{n_k} =u_k(k) \leq v_{n_k +1} \leq v_{n_k + 2} \leq \dots
\leq v_{n_k + k+1} \leq u_{k+1}(k) = u_{k+1}(k+1) = v_{n_{k+1}} .
$$
Set $u_i(t) = u_i(1) = v_{n_i}, t \in [k,k+1]$, when $i \geq k+1$.
Set $I_j = \left[ k + \frac{j}{k+1}, k + \frac{j+1}{k+1} \right]$,
$j = 0,1,2, \dots, k$.
On the interval $I_j$, $u_{k-j}(t),
t \in I_j$, is the straight line from $u_{k-j}(k)$ to $v_{n_k + k -j}$,
i.e.
$$
u_{k-j}(t) = (j+1 - (k+1)(t - k))u_{k-j}(k) +
((k+1)(t-k) - j)v_{n_k + k -j} ,
$$
$t \in I_j$ and other $u_m(t)$, $m\neq k-j$, are constants.
The construction of $\{u_i(t)\}_{i=0}^{\infty},
t \in [1,\infty)$, can then proceed by induction.

Set $\Delta_0(t)  = \sqrt{u_0(t)},$ $\Delta_i(t) = \sqrt{u_{i}(t) -
u_{i-1}(t)}$, $i \geq 1$. Let $\{\psi_{t_n}\}_{n=0}^{\infty}$ be a
discretization of $\psi$ such that $t_i \leq i$ for all $i \geq 1$. Set
$$
\Psi_t(a) = \sum_{i=0}^{\infty} \Delta_i(t) \psi_{\max \{t,t_i\} }(a)
\Delta_i(t).
$$
The sequence converges in the strict topology of $M(B)$ by
Lemma 3.1 of \cite{MT5}. Note that it follows from (\ref{4})
that for each $n \in \mathbb N$ there is an $N_n \in \mathbb N$
such that
\begin{equation*}
\begin{split}
&\Psi_t(x) - \Psi_s(x) =  \sum_{i = 0}^{N_n} \Delta_i(t)
\psi_{\max \{t,t_i\} }(x) \Delta_i(t) -
\Delta_i(s) \psi_{\max \{s,t_i\} }(x) \Delta_i(s)
\end{split}
\end{equation*}
for all $s,t \in [1,n]$. This shows that $\Psi_t(x) -
\Psi_s(x) \in B$ for all $s,t \in [1,\infty)$ and that
$t \mapsto \Psi_t(x)$ is continuous. We claim that
$\Psi_t(J) \subseteq B$ for all $t$. By Lemma 3.1 of \cite{MT5},
$\Psi_t, t \in [1,\infty)$, is an equi-continuous family since
$\psi_t, t \in [1,\infty)$, is, so it suffices to show that
$\Psi_t(x) \in B$ when $x \in F_k\cap J$. As we know that
$\Psi_t(x) - \Psi_1(x) \in B$, we must show that $\Psi_1(x) \in B$.
It follows from (\ref{3'}), (\ref{4}) and (\ref{2'}) that
$$
\|\Delta_i(1)\psi_{t_i}(x)\| \leq 2^{-i-1}
$$
when $i \geq k$. This shows that $\sum_{i=0}^{\infty} \Delta_i(1)
\psi_{t_i }(x) \Delta_i(1)$ converges in norm, proving that
$\Psi_1(x) \in B$.

To show that $\Psi$ is asymptotically multiplicative, it suffices,
by equi-continuity of $\Psi_t, t \in [1,\infty)$, to check for
$x,y \in F_k$. In the following we write $ a \sim_{\delta} b$ when
$a$ and $b$ are elements of the same $C^*$-algebra and
$\|a-b\| \leq \delta$. Let $t \in [m,m+1], \ m \geq k$.
When $i > m$, we have that $\max \{t,t_i\} \leq i$, and
hence
$$
\Delta_i(t)\psi_{\max\{t,t_i\}}(x)\Delta_i(t) \sim_{2^{-i}}
\psi_{\max\{t,t_i\}}(x)\Delta_i(t)^2,
$$
thanks to (\ref{4}), (\ref{3'}), (\ref{1}) and (\ref{1u}). Similarly,
$$
\Delta_i(t)\psi_{\max\{t,t_{i-1}\}}(x)\Delta_i(t) \sim_{2^{-i}}
\psi_{\max\{t,t_{i-1}\}}(x)\Delta_i(t)^2,
$$
and both estimates also hold with $x$ replaced by $y$. When
$i \leq m$, $\max\{t,t_i\} \leq m+1$, while
$u_i(t), u_{i-1}(t) \in \co \left\{ v_j : j \geq m \right\}$
by (\ref{3}). It follows therefore from (\ref{1}) and (\ref{1u}) that
$$
\Delta_i(t)\psi_{\max\{t,t_i\}}(x)\Delta_i(t) \sim_{2^{-m}}
\psi_{\max\{t,t_i\}}(x)\Delta_i(t)^2.
$$
Similarly,
$$
\Delta_i(t)\psi_{\max\{t,t_{i-1}\}}(x)\Delta_i(t)
\sim_{2^{-m}} \psi_{\max\{t,t_{i-1}\}}(x)\Delta_i(t)^2,
$$
and both estimates also hold with $x$ replaced by $y$. Set
$$
\delta_1(t) = \sup_j \|\psi_{\max\{t,t_j\}}(y) -
\psi_{\max\{t,t_{j+1}\}}(y) \|,
$$
$$
\delta_2(t) = \sup_j \|\psi_{\max\{t,t_j\}}(x)\psi_{\max\{t,t_j\}}(y) -
\psi_{\max\{t,t_j\}}(xy)\|,
$$
and
$$
k_x = \sup_t \|\psi_t(x)\|.
$$
Using Lemma 3.1 of \cite{MT5} and the above estimates we find that
 $$
\Psi_t(x)\Psi_t(y)\ \ = \ \
\left( \sum_{j=0}^{\infty} \Delta_j(t) \psi_{\max\{t,t_j\}}(x)
\Delta_j(t) \right) \left( \sum_{j=0}^{\infty} \Delta_j(t)
\psi_{\max\{t,t_j\}}(y)\Delta_j(t) \right)
\ \ \ \ \ \ \ \ \ \
 $$
\begin{eqnarray*}
&=&  \sum_{j=0}^{\infty}  \Delta_j(t) \psi_{\max\{t,t_j\}}(x)
\Delta_j(t)^2  \psi_{\max\{t,t_j\}}(y)\Delta_j(t)  \\
&& +\ \sum_{j=0}^{\infty}  \Delta_j(t)
\psi_{\max\{t,t_j\}}(x)  \Delta_j(t)\Delta_{j+1}(t)
\psi_{\max\{t,t_{j+1}\}}(y)\Delta_{j+1}(t) \\
&& +\  \sum_{j=0}^{\infty}  \Delta_{j+1}(t)
\psi_{\max\{t,t_{j+1}\}}(x)  \Delta_{j+1}(t)\Delta_{j}(t)
\psi_{\max\{t,t_{j}\}}(y)\Delta_{j}(t)
\end{eqnarray*}
\begin{eqnarray*}
& \sim_{6 k_xm2^{-m}} &  \sum_{j=0}^{\infty}  \Delta_j(t)
\psi_{\max\{t,t_j\}}(x)   \psi_{\max\{t,t_j\}}(y) \Delta_j(t)^2
\Delta_j(t)  \\
&& +\  \sum_{j=0}^{\infty}  \Delta_j(t)
\psi_{\max\{t,t_j\}}(x)  \psi_{\max\{t,t_{j+1}\}}(y)
\Delta_j(t)\Delta_{j+1}(t) \Delta_{j+1}(t) \\
&& +\  \sum_{j=0}^{\infty}  \Delta_{j+1}(t)
\psi_{\max\{t,t_{j+1}\}}(x)   \psi_{\max\{t,t_{j}\}}(y)
\Delta_{j+1}(t)\Delta_{j}(t) \Delta_{j}(t)
\end{eqnarray*}
\begin{eqnarray*}
& \sim_{2 k_x\delta_1(t)} &  \sum_{j=0}^{\infty}  \Delta_j(t)
\psi_{\max\{t,t_j\}}(x)   \psi_{\max\{t,t_j\}}(y) \Delta_j(t)^2
\Delta_j(t)  \\
&& +\  \sum_{j=0}^{\infty}  \Delta_j(t)
\psi_{\max\{t,t_j\}}(x)  \psi_{\max\{t,t_{j}\}}(y)
\Delta_j(t)\Delta_{j+1}(t) \Delta_{j+1}(t) \\
&& +\  \sum_{j=1}^{\infty}  \Delta_{j}(t)
\psi_{\max\{t,t_{j}\}}(x)   \psi_{\max\{t,t_{j}\}}(y)
\Delta_{j}(t)\Delta_{j-1}(t) \Delta_{j-1}(t)
\end{eqnarray*}
\begin{eqnarray*}
& \sim_{3\delta_2(t)}& \sum_{j=0}^{\infty}  \Delta_j(t)
\psi_{\max\{t,t_j\}}(xy)  \Delta_j(t)^2 \Delta_j(t)  \\
&& +\  \sum_{j=0}^{\infty}  \Delta_j(t)
\psi_{\max\{t,t_j\}}(xy)  \Delta_j(t)\Delta_{j+1}(t)
\Delta_{j+1}(t) \\
&& +\  \sum_{j=1}^{\infty}  \Delta_{j}(t)
\psi_{\max\{t,t_{j}\}}(xy)   \Delta_{j}(t)\Delta_{j-1}(t)
\Delta_{j-1}(t)
\end{eqnarray*}
 $$
 =\ \ \ \Psi_t(xy).\ \ \ \ \ \ \ \ \ \ \ \ \ \ \ \ \ \ \ \ \ \ \ \ \
\ \ \ \ \ \ \ \ \ \ \ \ \ \ \ \ \ \ \ \ \
 $$
Since $6k_x m2^{-m} + 2k_x\delta_1(t) + 3\delta_2(t)$ goes to zero as
$m$ tends to infinity, we conclude that $\lim_{t \to \infty}
\Psi_t(x)\Psi_t(y) - \Psi_t(xy) = 0$. Asymptotic linearity and
self-adjointness follow in the same way. Thus $\Psi$ is a
constantly extended asymptotic homomorphism. For $a \in A, s \in [0,1]$,
define $\Lambda_t(a)(s) \in M(B)$ by the strictly convergent sequence
$$
\Lambda_t(a)(s) = \sum_{i=0}^{\infty} \Delta_i(t)
\psi_{s \max \{t,t_i\} + (1-s)t }(a) \Delta_i(t).
$$
Since $s \mapsto \Lambda_t(a)(s)$ is a strictly continuous
and normbounded function, we have defined a family of maps
$\Lambda_t: A \to M(IB), t \in [1,\infty)$. It follows from
(\ref{4}) that for fixed $n$ there is an $N_n$ so large that
\begin{equation*}
\begin{split}
&\Lambda_t(a)(s) - \Lambda_{t'}(a)(s) = \\
&\sum_{i=0}^{N_n} \left(\Delta_i(t) \psi_{s \max \{t,t_i\} +
(1-s)t }(a) \Delta_i(t) -\Delta_i(t') \psi_{s \max \{t',t_i\} +
(1-s)t' }(a) \Delta_i(t')\right) \\
& \ \ \ \ \ \ \ \ \ \ \ \ \ +  \sum_{i=N_n+1}^{\infty} \Delta_i(1)
\left(\psi_{s \max \{t,t_i\} + (1-s)t }(a) - \psi_{s \max \{t',t_i\} +
(1-s)t' }(a)\right) \Delta_i(1)
\end{split}
\end{equation*}
for all $a$ and $s$, provided $t,t' \in [1,n]$. When $t$ tends to
$ t'$, the first term converges to $0$ in norm, uniformly in $s$,
for obvious reasons, and the second term does the same thanks to
Lemma 3.1 of \cite{MT5} and the continuity of $t \mapsto \psi_t(a)$.
Thus $t \mapsto \Lambda_t(a)$ is normcontinuous. Lemma 3.1 of
\cite{MT5} also shows that the family $\left(\Lambda_t\right)_{t
\in [1,\infty)}$, is equi-continuous since
$\left(\psi_t\right)_{t \in [1,\infty)}$ is. To show that
$\Lambda_t(J) \subseteq IB$ for all $t$, we must give an
argument different from the one used above since
$\Lambda_t - \Lambda_1$ does not map $J$ into $IB$.
Note that $s\max\{t,t_i\} + (1-s)t \leq \max\{t,i\} \leq i$
when $i \geq t$. According
to (\ref{4}) and (\ref{3'}), $u_i(t) = v_{n_i} \geq v_i$ when $i \geq t+1$, so we conclude from (\ref{2}) and
(\ref{2'}) that $\sup_s \|\Delta_i(t) \psi_{s \max \{t,t_i\} +
(1-s)t }(x)\| \leq 2^{-i}$ for all large enough $i$, when $x
\in \bigcup_k F_k \cap J$. Hence
the sum defining $\Lambda_t(x)$ converges in norm to an element of
$IB$. By continuity of $\Lambda_t$, we conclude that $\Lambda_t(J)
\subseteq IB$. The arguments that proved that $\Psi$ is an asymptotic
homomorphism show the same about $\Lambda$, thanks to the uniform
continuity of $\psi$. (The uniform continuity is used to show that the analog
of $\delta_1(t)$ tends to $0$ when $t$ goes to infinity.) 
$\Lambda$ is consequently an extended asymptotic
homomorphism given us a homotopy connecting $\Psi$ to
$\left(\sum_{i=0}^{\infty} \Delta_i(t)\psi_t(\cdot)
\Delta_i(t)\right)_{t \in [1,\infty)}$. For each
$a,\in A, s \in [0,1]$, define $\mu_t(a)(s) \in M(B)$ by
\begin{equation}\label{?2}
\mu_t(a)(s) = \begin{cases} \sum_{i=0}^{\infty} \Delta_i(t -
\log s)\psi_t(a) \Delta_i(t - \log s), &  s \neq 0, \\
\psi_t(a) , & s = 0. \end{cases}
\end{equation}
Since $\Delta_0(t)$ strictly tends to 1 and $\Delta_i(t)$, $i>0$, strictly
tend to 0, as $t\to\infty$,
the formula (\ref{?2}) defines an extended asymptotic
homomorphism $\mu = (\mu_t)_{t \in [1,\infty)} : A \to M(IB)$ providing a
homotopy between $\left(\sum_{i=0}^{\infty} \Delta_i(t)\psi_t(\cdot)
\Delta_i(t)\right)_{t \in [1,\infty)}$ and
$\left(\psi_t\right)_{t \in [1,\infty)}$.

Injectivity: Let $\varphi, \psi : A \to M(B)$ be constantly
extended asymptotic homomorphisms that are homotopic as
extended  asymptotic homomorphisms, and let $\Phi : A \to M(IB)$
be an extended asymptotic homomorphism realizing a homotopy
between the two. Disregarding a few considerations concerning
equi-continuity and uniform continuity, the construction from
the proof of surjectivity gives us a homotopy of constantly
extended asymptotic homomorphisms between
$\left( \sum_{j=0}^{\infty} \Delta_j(t)\varphi_t( \cdot)
\Delta_j(t)\right)_{t \in [1,\infty)}$ and $\left( \sum_{j=0}^{\infty}
\Delta_j(t)\psi_t( \cdot) \Delta_j(t)\right)_{t \in [1,\infty)}$,
where the $\Delta_i$'s arise from appropriately chosen continuous
paths of unit sequences in $B$. To complete the proof it suffices
therefore to check that the asymptotic homomorphism $\mu$ of
(\ref{?2}) is constantly extended when $\psi$ is. So assume
this is the case and let $\epsilon > 0$ be given. Let $t,t'
\in [1,\infty)$.
By equi-continuity it suffices to show that $\mu_t(x) - \mu_{t'}(x)
\in IB$ when $x \in F_k$. Take $n \in \mathbb N$ such that $n \geq
\max\{t,t'\}$. It follows then from (\ref{3}) that $u_i(t - \log s)
\in \co \{v_j : j \geq n\}$ for all $i \in \mathbb N$ and all
$s \in ]0,1]$. It follows therefore from (\ref{1}) and (\ref{1u})
that $\left\|\left[ \Delta_i(t - \log s), \psi_t(x)\right]\right\|
\leq 2^{-i}$ when $i \geq \max \{n,k\}$. As a consequence there is
an $N \in \mathbb N$ so large that
$$
 \psi_t(x)\left( 1 - u_{N}( t - \log s)\right) + \sum_{i=0}^{N}
\Delta_i(t - \log s)\psi_t(x) \Delta_i(t - \log s)  \
\sim_{\epsilon} \  \sum_{i=0}^{\infty} \Delta_i(t - \log s)
\psi_t(a) \Delta_i(t - \log s)
$$
for all $s \in ]0,1]$. By increasing $N$ we may assume that the
same estimate holds with $t$ replaced by $t'$. Thus $\mu_t(x) -
\mu_{t'}(x)$ has distance less than $2 \epsilon$ to the element
of $M(IB)$ given by the strictly continuous map $f: [0,1] \to M(B)$,
where
\begin{equation}\label{?4}
\begin{split}
& f(s) = \sum_{i=0}^{N} \left(\Delta_i(t - \log s)\psi_t(x)
\Delta_i(t - \log s) - \Delta_i(t' -  \log s)\psi_{t'}(x)
\Delta_i(t' - \log s)\right) \\
& \ \ \ \ \ \ \ \ \ \ \ \ + \psi_t(x)\left( 1 - u_{N}( t -
\log s)\right) - \psi_{t'}(x)\left( 1 - u_{N}( t' - \log s)\right) ,
\end{split}
\end{equation}
when $s \in ]0,1]$, and
\begin{equation}\label{?5}
f(s) = \psi_t(x) - \psi_{t'}(x) ,
\end{equation}
when $s = 0$. Note that for each $s$, $f(s)$ is in $B$ since
$\psi_t(x) - \psi_{t'}(x)$ is, and that $f$ is obviously
norm-continuous on $]0,1]$. It suffices now to show that
(\ref{?4}) converges in norm to (\ref{?5}) when $s$ tends to
zero. To see that this is the case note that (\ref{3}),
(\ref{1}) and (\ref{1u}) imply that
$$
\lim_{s \to 0} \left( \sum_{i=0}^{N} \Delta_i(t - \log s)\psi_t(x)
\Delta_i(t - \log s) -  \sum_{i=0}^{N} \psi_t(x) \left(
\Delta_i(t - \log s)\right)^2\right) = 0 .
$$
The same conclusion holds with $t$ replaced by $t'$ so (\ref{?4})
approaches (\ref{?5}) as $s \to 0$ because $ \sum_{i=0}^{N}
\left(\Delta_i(t - \log s)\right)^2 + \left( 1 - u_{N}( t -
\log s)\right) =  \sum_{i=0}^{N}  \left(\Delta_i(t' - \log s)
\right)^2 + \left( 1 - u_{N}( t' - \log s)\right) = 1$ for all
$s \in ]0,1]$. The proof is complete.
\end{proof} \qed

Theorem \ref{Thm3} serves as our excuse for not distinguishing
very strictly between $[[A,J;B]]$ and $[\{A,J;B\}]$ in the following.

When $B$ is stable both $[[A,J;B]]$ and $[\{A,J;B\}]$ are equipped
with a semi-group structure in the familiar way: When $\varphi, \psi :
A \to M(B)$ are (constantly) extended asymptotic homomorphisms and
$V_1,V_2 \in M(B)$ are isometries such that $V_1V_1^* + V_2V_2^* = 1$,
we can define a (constantly) extended asymptotic homomorphism $\varphi
\oplus \psi : A \to M(B)$ by $(\varphi \oplus \psi)_t(a) =
V_1\varphi_t(a)V_1^* + V_2\psi_t(a)V_2^*$. The compositions
defined in this way in $[[A,J;B]]$ and $[\{A,J;B\}]$ are
commutative and associative, and are independent of the choice
of isometries $V_1,V_2$. In this case the bijection of Theorem
\ref{Thm3} is an isomorphism of abelian semi-groups. In the
following we assume that $B$ is stable.

\begin{lem}\label{triv} $[[A,0;B]] = 0$.
\end{lem}
\begin{proof} Consider an asymptotic homomorphism $\pi : A \to M(B)$
such that $\pi_t(0) \in B$ for all $t$. Then $\pi'_t(a) =
\pi_t(a) - \pi_t(0)$ defines an asymptotic homomorphism $\pi' :
A \to M(B)$ with the property that $\pi'_t(0) = 0$ for all $t$,
and $[\pi] = [\pi']$ in $[[A,0;B]]$. By Lemma 1.3.6 of \cite{K-JT}
there is a strictly continuous family $V_s, s \in ]0,1]$, of
isometries in $M(B)$ such that $V_1 =1$ and
$\lim_{s \to 0} V_sV_s^* = 0$ in the strict topology. Set
$$
\Phi_t(a)(s) = \begin{cases} V_s\pi'_t(a)V_s^* , & s \in ]0,1], \\
0, & s = 0. \end{cases}
$$
Note that $s \mapsto \Phi_t(a)(s)$ is strictly continuous and
norm-bounded. Thus $\Phi = \left(  \Phi_t\right)_{t \in [1,\infty)} :
A \to M(IB)$ is an asymptotic homomorphism such that $\Phi_t(0) = 0$
for all $t$, giving us a homotopy connecting $\pi'$ to $0$.
\end{proof} \qed

We denote the $C^*$-algebras $C(\mathbb T)\otimes A$ and
$C_0(0,1) \otimes A$ by $TA$ and $SA$, respectively. Note that
there is an extension
\begin{equation}\label{?1}
\begin{xymatrix}{
0 \ar[r] & SA \ar[r] & TA \ar[r]^{\ev} & A \ar[r] & 0 , }
\end{xymatrix}
\end{equation}
where $\ev : TA \to A$ is evaluation at $1 \in \mathbb T$. We shall
often identify $TA$ with $\{ f \in IA : f(0) = f(1)\}$ in the obvious way.

\begin{lem}\label{?1?}
$[[TA,SA;B]]$ and $[\{TA,SA;B\}]$ are groups.
\end{lem}
\begin{proof} Since the bijection of Theorem \ref{Thm3} is an
isomorphism of semi-groups with zero, it suffices to show that
$[[TA,SA;B]]$ is a group. Let $\varphi: TA \to M(B)$ be an
extended asymptotic homomorphism. Let $\alpha \in \Aut TA$
be the automorphism which changes orientation on the circle, i.e.
$\alpha (f)(s) = f(1-s)$, $f\in T$, $s\in[0,1]$.
Then the homomorphism $\gamma : TA \to M_2(TA)$ given
by
$$
\gamma(f) =  \left( \begin{smallmatrix} f  &  \\ & \alpha(f)
\end{smallmatrix} \right)
$$
is homotopic to the $*$-homomorphism
$$
f \mapsto   \left( \begin{smallmatrix} f(0)  &  \\ & f(0)
\end{smallmatrix} \right),
$$
via a path of $*$-homomorphisms which all send $SA$ into $M_2(SA)$.
Since $ \varphi \oplus (\varphi \circ \alpha) = \tilde{\Theta} \circ \left(
\id_{M_2} \otimes \varphi\right)  \circ \gamma$, where the
$*$-isomorphism $\tilde{\Theta} :M_2(M(B)) \to M(B)$ is given by
\begin{equation}\label{stab}
\tilde{\Theta} \left( \begin{smallmatrix} a_{11} & a_{12} \\ a_{21} & a_{22}
\end{smallmatrix} \right) = V_1a_{11}V_1^* + V_1a_{12}V_2^* +
V_2a_{21}V_1^* + V_2 a_{22}V_2^*,
\end{equation}
we conclude that $ \varphi \oplus (\varphi \circ \alpha)$ is homotopic
as an extended asymptotic homomorphism to $\left(\varphi \circ c
\circ \ev\right) \oplus \left(\varphi \circ c \circ \ev\right)$,
where the $*$-homomorphisms $c: A \to TA$ and $\ev:TA \to A$ are
given by $c(a)(t) = a, t \in \mathbb T$, and $\ev(f) = f(0)$,
respectively. Since $[\varphi \circ c] = 0$ in $[[A,0;B]]$ by Lemma
\ref{triv}, it follows that $[\left(\varphi \circ c
\circ \ev\right) \oplus \left(\varphi \circ c \circ \ev\right)] = 0$ in
$[[TA,SA;B]]$.
\end{proof} \qed

\section{Various maps}

Let $A$ and $B$ be separable $C^*$-algebras, $B$ stable. In this
section we obtain our main result which is that an appropriate
modification of the Connes-Higson construction gives rise to an
isomorphism between $\Ext^{-1/2}(A,B)$ and $[[TA,SA; B]]$.

\subsection{The Connes--Higson map}\label{CHmap}
 Let $\phi:A\to Q(B)$ be a semi-invertible extension.
Then there exists
an extension $\psi:A\to Q(B)$ and an equi-continuous asymptotic homomorphism
$\pi=(\pi_t)_{t\in[1,\infty)}:A\to
M_2(M(B))$ such that $q_{M_2(B)}\circ\pi=\phi\oplus\psi$. Denote the
matrix
elements of $\pi_t$ by $\pi_t^{ij}$, $i,j=1,2$. Note that
$\pi^{12}_t(A) \cup \pi^{21}_t(A) \subseteq   B$ for all $t$. It
follows from the equi-continuity of
$\pi$ and the separability of $A$ and $B$ that there exists an
 approximate unit $(u_t)_{t\in[1,\infty)}\subseteq B$
such that
 \begin{equation}\label{CH1}
\lim_{t\to\infty}\left[f(u_t),\pi_t^{11}(a)\right]=0,
\end{equation}
\begin{equation}\label{CH2}
\lim_{t\to\infty}\pi_t^{12}(a) \left(f(u_t)-f(1)\right)=0,
\end{equation}
and
\begin{equation}\label{CH3}
\lim_{t\to\infty}\pi_t^{21}(a)\left(f(u_t)-f(1)\right)=0,
\end{equation}
for all $f \in C[0,1]$ and all $a \in A$.
Then $\left( \begin{smallmatrix} f(u_t) &  \\ & f(0) \end{smallmatrix}
\right)$ and $\pi_t(a)$ asymptotically commute for all $f \in T$ and
$a \in A$. We use here and in the following $T$ to denote the
$C^*$-algebra $C(\mathbb T)$. Note that
$$
\pi_t(a) - \pi_1(a) \in M_2(B)
$$
and
$$
\left( \begin{smallmatrix} f(u_t) &  \\ & f(0) \end{smallmatrix} \right)
- \left( \begin{smallmatrix} f(0) &  \\ & f(0) \end{smallmatrix} \right)
\in M_2(B)
$$
for all $a \in A, \ f \in T$. Set
$$
X = \{g \in C_b([1,\infty),M_2(M(B))): g(t)-g(1) \in M_2(B) \ \forall t \}.
$$
It follows that there is a $*$-homomorphism $\Phi : TA \to
X/C_0([1,\infty),M_2(B))$ such that $\Phi(f \otimes a)$ is the
image of the element in $X$ given by the function
$$
t \mapsto \left( \begin{smallmatrix} f(u_t) &  \\ & f(0)
\end{smallmatrix} \right)\pi_t(a).
$$
It follows from the Bartle-Graves selection theorem that there is a
continuous
map $\chi : X/C_0([1,\infty),M_2(B)) \to X$ which is a right-inverse
for the quotient map $X \to  X/C_0([1,\infty),M_2(B)))$. By Remark
2 on page 114 of \cite{L} we can assume that $\chi$ maps the asymptotic
algebra of $M_2(B)$, which is a $C^*$-subalgebra of $X$, into
$C_b\left( [1,\infty), M_2(B)\right)$. Set
$$
CH(\varphi)_t(a) = \tilde{\Theta} \left( \left( \chi \circ
\Phi(a)(t)\right)\right),
$$
where $\tilde{\Theta} : M_2(M(B)) \to M(B)$ is the $*$-isomorphism
(\ref{stab}). Note that
$CH(\varphi)$ is an equi-continuous asymptotic homomorphism
$CH(\varphi) : TA \to M(B)$ such that
$CH(\varphi)_t(SA) \subseteq B$ and $CH(\varphi)_t(x) - CH(\varphi)_1(x)
\in B$
for all $t$ and all $x \in TA$. In short, $CH(\varphi)$ is an asymptotic
homomorphism which is constantly extended from $SA$,
and defines an element of $[\{TA,SA;B\}]$. It is easy to see that the
construction gives us a well-defined group homomorphism
$$
CH : \Ext^{-1/2}(A,B) \to [\{TA,SA;B\}].
$$
When composed with the obvious forgetful map $[\{TA,SA;B\}] \to [[SA,B]]$
obtained by restricting asymptotic homomorphisms to $SA$, we get the
usual Connes-Higson map.

\subsection{The $E$-map}\label{subsection_E}

Let $\varphi=(\varphi_t)_{t\in[1,\infty)} : TA \to M(B)$ be an
asymptotic homomorphism which is constantly extended from $SA$,
i.e. $\varphi$ is an asymptotic homomorphism such that $\varphi_t(SA)
\subseteq B$ and $\varphi_t(x) - \varphi_1(x) \in B$ for all
$t \in [1,\infty)$ and all
$x \in TA$. By Lemma \ref{gen2} we may assume that $\varphi$
is equi-continuous. We will use $\varphi$ to define a
semi-invertible extension of $T^2A$ by $B$, where $T^2A = T(TA) =
C\left(\mathbb T^2\right) \otimes A$. To do this we choose first
a discretization $\varphi_{t_0}, \varphi_{t_1}, \varphi_{t_2},
\dots$ such that $\lim_{i \to \infty} t_i = \infty $ and
$\lim_{i \to \infty} \sup_{t  \in [t_i,t_{i+1}]}
\|\varphi_t(a) - \varphi_{t_i}(a)\| = 0$ for all $a \in TA$.
To define from such a discretization a map $\bold\Phi :
TA \to  M(B \otimes \mathbb K)$ we identify $\mathbb K$ with
the compact operators on the Hilbert space $l^2(\mathbb Z)$,
and introduce the corresponding matrix units $e_{i,j} \in
\mathbb K, i,j \in \mathbb Z$. Then both sums in
\begin{equation}\label{klaus1}
\bold\Phi (f) = \sum_{i \geq 1} \varphi_{t_i}(f) \otimes e_{i,i} +
\sum_{i \leq 0} \varphi_{t_{|i|}}(f(0)) \otimes e_{i,i}
\end{equation}
converge in the strict topology and (\ref{klaus1}) defines a map $\bold\Phi :
TA \to M(B \otimes \mathbb K)$. Observe that $\bold\Phi$ is a
$*$-homomorphism modulo $B\otimes \mathbb K$. Furthermore,
$\bold\Phi(a)$ commutes modulo $B\otimes \mathbb K$ with the
two-sided shift $\mathcal T = \sum_{j \in \mathbb Z} e_{j,j-1}$.
So we get in this way an extension
$$
E(\varphi) : T^2A \to Q(B\otimes \mathbb K)
$$
such that
\begin{equation}\label{klaus2}
E(\varphi)(g \otimes f) = q_{B \otimes \mathbb K}
\left(g\left(\mathcal T\right)\bold\Phi(f)\right)
\end{equation}
for all $g \in T, f \in TA$.

\begin{lem}\label{LL2} $E( \varphi)$ is semi-invertible, and
its class in $\Ext^{-1/2}(T^2A,B)$ does not depend on the chosen
discretization of $\varphi$.
\end{lem}
\begin{proof} The inverse $-E(\varphi)$ is given by the formula
\begin{equation}\label{klaus3}
-E(\varphi)(g \otimes f) = q_{B \otimes \mathbb K}\left(g(\mathcal T)
\bold\Psi(f)\right)
\end{equation}
where
\begin{equation}\label{klaus4}
\bold\Psi(f) = \sum_{i \leq 0} \varphi_{t_{|i|}}(f)\otimes e_{i,i} +
\sum_{i \geq 1} \varphi_{t_i}(f(0)) \otimes e_{i,i}.
\end{equation}
To see that $\hat{\Theta}^{-1} \circ \left(E(\varphi) \oplus
\left( - E(\varphi)\right)\right)$ is asymptotically split it
is appropriate to view $M_2(B \otimes \mathbb K)$ as $M_2(B)\otimes
\mathbb K$. Let $z$ denote the identity function on the circle
$\mathbb T$, so that $z$ generates $T$ as a $C^*$-algebra. Then
$\hat{\Theta}^{-1} \circ \left( E(\varphi) \oplus
\left( - E(\varphi)\right)\right)$ is determined by the condition
that $\hat{\Theta}^{-1} \circ \left(E(\varphi) \oplus
\left( - E(\varphi)\right)\right) =q_{M_2(B \otimes \mathbb K)}
\circ \Psi$, where
\begin{equation*}
\begin{split}
&\Psi(z^k \otimes f) = \\
&\left( \sum_{i \in \mathbb Z} \left( \begin{smallmatrix} 1 & 0 \\
0 & 1 \end{smallmatrix} \right) \otimes e_{n,n-1} \right)^k
\left( \sum_{i  \geq 1}\left(  \begin{smallmatrix} \varphi_{t_i}(f) & 0 \\
0 & \varphi_{t_i}(f(0)) \end{smallmatrix} \right) \otimes e_{i,i} +
\sum_{i  \leq 0}  \left( \begin{smallmatrix}
\varphi_{t_{|i|}}(f(0)) & 0 \\ 0 & \varphi_{t_{|i|}}(f)
\end{smallmatrix} \right) \otimes e_{i,i}\right),
\end{split}
\end{equation*}
modulo $B \otimes \mathbb K$ for all $k \in \mathbb Z$ and all
$f \in TA$. Set
$$
\varphi^i_t = \begin{cases} \varphi_{\max\{t_i,t\}}, & i \geq 0, \\
\varphi_{\max\{t_{|i|},t\}}, & i \leq 0. \end{cases}
$$
Without loss of generality we may assume that the discretization
$\{t_i\}_{i=0}^\infty$ satisfies $t_0=t_1=1$, so that
$\varphi^0_t=\varphi^1_t=\varphi_t$.
Define a continuous family of unitaries, $S_t$, by
$$
S_t =  \sum_{n \neq 1} \left( \begin{smallmatrix} 1 & 0 \\
0 & 1 \end{smallmatrix} \right) \otimes e_{n,n-1}  +
\left( \begin{smallmatrix} 1 - u_t & \sqrt{2u_t - u_t^2} \\
\sqrt{2u_t - u_t^2} & 1 - u_t \end{smallmatrix} \right) \otimes e_{1,0} ,
$$
where $u_t, t \in [1,\infty)$, is a continuous approximate unit in $B$. Then
 $$
\left[S_t, \left( \sum_{i  \geq 1}\left(  \begin{smallmatrix}
\varphi^i_{t}(f) & 0 \\ 0 & \varphi^i_{t}(f(0)) \end{smallmatrix}
\right) \otimes e_{i,i} + \sum_{i  \leq 0}  \left( \begin{smallmatrix}
\varphi^i_{t}(f(0)) & 0 \\ 0 & \varphi^i_{t}(f) \end{smallmatrix}
\right) \otimes e_{i,i}\right)\right]
\ \ \ \ \ \ \ \ \ \ \ \ \ \ \ \ \ \ \ \ \ \
\ \ \ \ \ \ \ \ \ \ \ \ \ \ \ \ \ \ \ \ \ \
 $$
\begin{eqnarray}\label{klaus6}
&=&
\sum_{i\leq 0}\left( \begin{matrix} \varphi^{i-1}_{t}(f) -
\varphi^i_{t}(f) & {0} \\ {0} & \varphi^{i-1}_{t}(f(0)) -
\varphi^i_{t}(f(0)) \end{matrix} \right) \otimes e_{i,i-1}
\nonumber\\
&+&
\sum_{i>1}\left( \begin{matrix} \varphi^{i-1}_{t}(f(0)) -
\varphi^i_{t}(f(0)) &
{0} \\ {0} & \varphi^{i-1}_{t}(f) - \varphi^i_{t}(f) \end{matrix} \right)
\otimes e_{i,i-1}
\\
&+& \left( \begin{matrix}
(1-u_t)\varphi_{t}(f)-\varphi_{t}(f(0))(1-u_t)&
{[}\sqrt{2u_t - u_t^2},\varphi_{t}(f(0)){]} \\
{[}\sqrt{2u_t - u_t^2},\varphi_{t}(f){]} &
(1-u_t)\varphi_{t}(f(0))-\varphi_{t}(f)(1-u_t) \end{matrix} \right)
\otimes e_{1,0}. \nonumber
\end{eqnarray}
The first two terms in the right-hand side of (\ref{klaus6}) vanishes as
$t\to\infty$ due to the choice of the discretization.
Since $A$ is separable and $\varphi_t(f) - \varphi_t(f(0)) \in B$,
we can choose $\left( u_t\right)_{t \in [1,\infty)}$ such that
$$
\lim_{t \to \infty} \left[u_t,\varphi_t(f)\right] = \lim_{t \to \infty}
\left[u_t,\varphi_t(f(0))\right] = 0
$$
and
$$
\lim_{t \to \infty} \left(1 -u_t\right)\left(\varphi_t(f) -
\varphi_t(f(0))\right) =  0
$$
for all $f \in TA$. Such a choice ensures that the last term in
(\ref{klaus6}) also vanishes as $t\to\infty$, so we see that $S_t$
asymptotically commutes with $ \sum_{i  \geq 1}\left(
\begin{smallmatrix} \varphi^i_{t}(f) & 0 \\ 0 & \varphi^i_{t}(f(0))
\end{smallmatrix} \right) \otimes e_{i,i} + \sum_{i \leq 0}
\left( \begin{smallmatrix}  \varphi^i_{t}(f(0)) & 0 \\ 0 &
\varphi^i_{t}(f) \end{smallmatrix} \right) \otimes e_{i,i}$.
The last expression, as well as $S_t$, is constant in $t$ when
taken modulo $B \otimes \mathbb K$, so we can define an asymptotic
splitting $\Lambda$ for $\hat{\Theta}^{-1} \circ \left(E(\varphi)
\oplus \left( - E(\varphi)\right)\right)$ such that
$\Lambda_t(z^k \otimes f)$ asymptotically agrees with
$$
S_t^k\left( \sum_{i  \geq 1}\left(  \begin{smallmatrix}
\varphi^i_{t}(f) & 0 \\ 0 & \varphi^i_{t}(f(0)) \end{smallmatrix}
\right) \otimes e_{i,i} + \sum_{i  \leq 0}  \left( \begin{smallmatrix}
\varphi^i_{t}(f(0)) & 0 \\ 0 & \varphi^i_{t}(f) \end{smallmatrix}
\right) \otimes e_{i,i}\right)
$$
for all $k \in \mathbb Z$ and all $f \in TA$. The method used to
reach this conclusion will be used several times in the following,
so we give a detailed account here: Set
$$
X = \{ f \in C_{b}\left( [1,\infty), M_2(M(B \otimes \mathbb K))
\right) : \ f(1) - f(t) \in M_2(B \otimes \mathbb K) \ \forall t \},
$$
which is a $C^*$-algebra containing $C_0\left( [1,\infty),
M_2(B \otimes \mathbb K)\right)$ as an ideal. Since $S_t$
asymptotically commutes with $ \sum_{i  \geq 1}\left(
\begin{smallmatrix} \varphi^i_{t}(f) & 0 \\ 0 & \varphi^i_{t}(f(0))
\end{smallmatrix} \right) \otimes e_{i,i} + \sum_{i  \leq 0}
\left( \begin{smallmatrix}  \varphi^i_{t}(f(0)) & 0 \\ 0 &
\varphi^i_{t}(f) \end{smallmatrix} \right) \otimes e_{i,i}$,
which is an asymptotic homomorphism, we get straightforwardly a
$*$-homomorphism $\Phi$ from $T^2A$ into the asymptotic algebra
of $M_2(M(B\otimes \mathbb K))$ such that
\begin{equation}\label{ups}
\Phi(z^k \otimes f) = S_t^k   \left( \sum_{i  \geq 1}
\left(  \begin{smallmatrix} \varphi^i_{t}(f) & 0 \\ 0 &
\varphi^i_{t}(f(0)) \end{smallmatrix} \right) \otimes e_{i,i} +
\sum_{i  \leq 0}  \left( \begin{smallmatrix}  \varphi^i_{t}(f(0)) &
0 \\ 0 & \varphi^i_{t}(f) \end{smallmatrix} \right) \otimes e_{i,i}\right),
\end{equation}
modulo $C_0([1,\infty), M_2(M(B \otimes \mathbb K)))$, for all
$k \in \mathbb Z, f \in TA$. Since the right-hand side of (\ref{ups})
is constant in $t$, modulo $M_2(B \otimes \mathbb K)$, it follows
that $\Phi$ takes values in $X/C_0([1,\infty),
M_2(B \otimes \mathbb K))$, which is - or should be considered as -
a $C^*$-subalgebra of the asymptotic algebra of
$M_2(M(B \otimes \mathbb K))$. By the Bartle-Graves selection
theorem there is a continuous section $\chi : X/C_0([1,\infty),
M_2(B \otimes \mathbb K)) \to X$ for the quotient map
$X \to X/C_0([1,\infty), M_2(B \otimes \mathbb K))$.
Set $\Pi_t(x) = \chi \circ \Phi(x)(t)$. Then $\Pi$ is
an asymptotic homomorphism such that $q_{M_2(B \otimes \mathbb K)}
\circ \Pi_t = \hat{\Theta}^{-1} \circ \left(E(\varphi)
\oplus (-E(\varphi))\right)$ for all $t$.

It follows that $E(\varphi)$ is semi-invertible. That its class
in $\Ext^{-1/2}(T^2A,B)$ is independent of the choice of
discretization follows from the homotopy invariance of
$\Ext^{-1/2}$, Theorem \ref{Thm1}, by using, for example,
the construction of homotopy from Lemma 5.3 in \cite{MT1}.
\end{proof} \qed

It follows from Lemma \ref{LL2} that there is a well-defined
group homomorphism
$$
E : [\{TA,SA;,B\}] \to \Ext^{-1/2}(T^2A,B)
$$
such that $E[\varphi] = [E(\varphi)]$.

\begin{rmk}\label{remark} \rm

For use in arguments below we give another proof of the
semi-invertibility of $E(\varphi)$, i.e. of the fact that
that $\hat{\Theta}^{-1} \circ \left( E(\varphi) \oplus
(-E(\varphi))\right)$ is asymptotically split. For each
$m \in \mathbb N$, define a sequence $\alpha_i(m), i \in \mathbb Z$,
of real numbers as follows:
$$
\alpha_i(m) = \begin{cases} 1, & \ i \geq 0 \\ 1 + \frac{i}{m} , &
\ 0 > i > - m \\ 0 , & \ i \leq  -m \end{cases}
$$
For $t \in [m,m+1[$, set
\begin{equation*}
\alpha_i(t) =   (m+1 - t)\alpha_i(m)  +  (t -m)\alpha_i(m +1).
\end{equation*}
For each $t \in [1,\infty)$, set
$$
U_{t} = \sum_{i \in \mathbb Z} \left( \begin{matrix}
\sqrt{\alpha_i(t)} & - \sqrt{1 - \alpha_i(t)} \\
\sqrt{1 - \alpha_i(t)} & \sqrt{\alpha_i(t)} \end{matrix} \right)
\otimes e_{i,i} .
$$
Each $U_{t}$ is a unitary in $M_2(M(B \otimes \mathbb K))$ and
\begin{equation}\label{ogg8}
\lim_{t \to \infty} \left[ U_{t},  \sum_{i \in \mathbb Z}
\left( \begin{smallmatrix} 1 & 0 \\ 0 & 1 \end{smallmatrix} \right)
\otimes e_{i,i-1} \right] = 0.
\end{equation}
Set
$$
\Gamma_t(f)   = \left( \sum_{i  \geq 1}\left(  \begin{smallmatrix}
\varphi^i_{t}(f) & 0 \\ 0 & \varphi^i_{t}(f(0)) \end{smallmatrix}
\right) \otimes e_{i,i} + \sum_{i  \leq 0}  \left( \begin{smallmatrix}
\varphi^i_{t}(f) & 0 \\ 0 & \varphi^i_{t}(f(0)) \end{smallmatrix}
\right) \otimes e_{i,i}\right)
$$
and note that
$$
U_{t}\Gamma_t(f)U_{t}^* -  \left(\sum_{i  \geq 1}
\left(  \begin{smallmatrix} \varphi_{t_i}(f) & 0 \\
0 & \varphi_{t_i}(f(0)) \end{smallmatrix} \right) \otimes e_{i,i} +
\sum_{i  \leq 0}  \left( \begin{smallmatrix}
\varphi_{t_{|i|}}(f(0)) & 0 \\ 0 & \varphi_{t_{|i|}}(f)
\end{smallmatrix} \right) \otimes e_{i,i} \right)
\in M_2(B \otimes \mathbb K) .
$$
Then
\begin{equation}\label{rem1}
\hat{\Theta}^{-1} \circ \left(E(\varphi)
\oplus (-E(\varphi))\right) (z^k \otimes f) =
q_{M_2(B \otimes \mathbb K)} \left( \left(\sum_{i \in \mathbb Z}
\left( \begin{smallmatrix} 1 & 0 \\ 0 & 1 \end{smallmatrix} \right)
\otimes e_{n,n-1}  \right)^k U_{t}\Gamma_t(f)U_{t}^*\right),
\end{equation}
for all $k \in \mathbb Z, f \in TA$ and all $t$. By the method
used at the end of the proof of Lemma \ref{LL2} we see from this
that $\hat{\Theta}^{-1} \circ \left(E(\varphi) \oplus
(-E(\varphi))\right)$ is asymptotically split.
\end{rmk}

\subsection{The Bott maps}\label{Bott-section}

We need a version of the Bott isomorphism in one of its many guises.
The one which best serves
our purpose is based on a particular projection in $M_2(T^2)$
which we now describe.

Given two commuting unitaries $S,T$ in a $C^*$-algebra, we define a
projection $P(S,T)$ in the $2 \times 2$ matrices over the
$C^*$-algebra generated by $S$ and $T$ in the following way.
Let $s, c_0,c_1 : [0,1] \to \mathbb R$ be the functions
$$
c_0(t) = |\cos (\pi t)|1_{[0,\frac{1}{2}]}(t), \ c_1(t) =
|\cos (\pi t)|1_{(\frac{1}{2},1]}(t), \ s(t) = \sin (\pi t) .
$$
Set $\tilde{g} = sc_0, \ \tilde{h} = sc_1$ and $\tilde{f} = s^2$.
Since $\tilde{f}, \ \tilde{g}$ and $\tilde{h}$ are continuous and
$1$-periodic they give rise to continuous functions, $f,g,h$, on
$\mathbb T$ (we identify $\tilde{f}$ on $[0,1]$ with $f$ on $\mathbb T$
in such a way that if $S=e^{2\pi ix}$ then $f(S)=\tilde{f}(x)$).
Set
$$
P(S,T)  =  \left ( \begin{matrix} f(S) & g(S) + h(S)T \\
T^*h(S) + g(S) & 1 - f(S) \end{matrix} \right ),
$$
cf. \cite{L}. When we apply the recipe to the canonical
generating unitaries of $C(\mathbb T^2)$, we get the desired
projection $P \in C(\mathbb T^2)\otimes M_2$. Then $P$ takes
the form
\begin{equation}\label{Bproj}
P(t,z) =   \left ( \begin{matrix} f(t) & g(t) + h(t)z \\
h(t)\overline{z} + g(t) & 1 - f(t) \end{matrix} \right ),
\end{equation}
$t \in [0,1], z \in \mathbb T$. Let $P_0 =
\left( \begin{smallmatrix} 0 & 0 \\ 0 & 1 \end{smallmatrix}
\right) \in M_2(T^2)$. Given a semi-invertible extension
$\phi : T^2A \to Q(B)$ we set
\begin{equation}\label{01}
\phi_1(a)=\hat{\Theta}\left((\id_{M_2}\otimes\phi)(P\otimes a)\right),
\qquad
\phi_0(a)=\hat{\Theta}\left((\id_{M_2}\otimes\phi)(P_0\otimes a)\right),
 \end{equation}
$a\in A$, where $\hat{\Theta}: M_2(Q(B)) \to Q(B)$ is the
isomorphism induced by the isomorphism $\tilde{\Theta}$
(\ref{stab}). It is
easy to see that $\varphi_1$ and $\varphi_0$ are both
semi-invertible since $\phi$ is. We define the Bott map
$\Bott : \Ext^{-1/2}(T^2A,B) \to \Ext^{-1/2}(A,B)$ such that
 $$
\Bott(\phi)=[\phi_1]-[\phi_0],
 $$
where $[\phi_i]$, $i=0,1$, are the classes of $\phi_i$ in
$\Ext^{-1/2}(A,B)$.

We can also use the projections $P$ and $P_0$ to define a Bott map
$$
\Bott : [[T^3A,ST^2A;B]] \to [[TA,SA;B]]
$$
such that $\Bott([\varphi]) = [\varphi^1] - [\varphi^0]$, where
$$
\varphi^1_t(a) = \tilde{\Theta}\left((\id_{M_2} \otimes \varphi_t)
\left(P \otimes a\right)\right),
$$
and
$$
\varphi^0_t(a) = \tilde{\Theta}\left((\id_{M_2} \otimes \varphi_t)
\left(P_0 \otimes a\right)\right)
$$
for all $t$ and all $a \in TA$. It is easy to see that the diagram
\begin{equation}\label{natdiag}
\begin{xymatrix}{
\Ext^{-1/2}(T^2A, B) \ar[d]_-{CH} \ar[r]^-{\Bott} & \Ext^{-1/2}(A,B)
\ar[d]^-{CH} \\
[[T^3A,ST^2A;B]] \ar[r]^-{\Bott} & [[TA,SA; B]] }
\end{xymatrix}
\end{equation}
commutes.

\section{From semi-invertible extensions to asymptotic homomorphisms
and back}\label{mainthm}
 In this section we prove our main result which is the following theorem.

\begin{thm}\label{mainth} Let $A$ and $D$ be separable $C^*$-algebras. Then
$$
CH : \Ext^{-1/2}(A,D) \to [[TA,SA;D \otimes \mathbb K]]
$$
is an isomorphism.
\end{thm}

We will prove Theorem \ref{mainth} by establishing the commutativity of the following
two diagrams:
 \begin{equation}\label{inj}
\begin{xymatrix}{%
\Ext^{-1/2}(A,B)\ar[rr]^-{CH}\ar@/_1pc/[drrrr]_-{e}&&
[[TA,SA;B]]\ar[rr]^-{E}&&
\Ext^{-1/2}(T^2A,B\otimes \mathbb K)\ar[d]^{Bott}\\
&&&&\Ext^{-1/2}(A,B \otimes \mathbb K),
}
\end{xymatrix}
\end{equation}
and
\begin{equation}\label{surj}
\begin{xymatrix}{%
[[TA,SA;B]]\ar[rr]^-{E}\ar@/_1pc/[drrrr]_-{e}&&
\Ext^{-1/2}(T^2A,B\otimes\mathbb K)\ar[rr]^-{Bott}&&
\Ext^{-1/2}(A,B\otimes\mathbb K)\ar[d]^-{CH}\\
&&&&[[TA,SA;B\otimes\mathbb K]],
}
\end{xymatrix}
\end{equation}
where $e $ in both cases is an isomorphism induced by the
stabilizing map $b \mapsto b \otimes e_{11}$ for some minimal
non-zero projection $e_{11} \in \mathbb K$. From the commutativity of
the first diagram we conclude that $CH$ is injectivity, and from the 
commutativity of the latter that $CH$ is surjective.

\subsection{$CH$ is injective}\label{iso}

For simplicity of notation we shall ignore the $*$-isomorphism
$\Theta$ in the definition of $CH$, and consider instead $CH$
as a map $CH: \Ext^{-1/2}(A,B) \to [[TA,SA;M_2(B)]]$.
Similarly, we will consider $\Bott$ as a map $\Bott :
\Ext^{-1/2}(T^2A,B) \to \Ext^{-1/2}(A,M_2(B))$. Let
$\varphi \in \Hom(A,Q(B))$ be a semi-invertible extension.
There is then an equi-continuous and essentially constant
asymptotic homomorphism
$$
\left( \begin{smallmatrix} \alpha_t & \beta_t \\ \gamma_t &
\delta_t \end{smallmatrix} \right)_{t \in [1,\infty)} :
A \to M_2\left(M(B)\right)
$$
such that $\varphi = q_B \circ \alpha_t$ for all $t$. By
'essentially constant' we refer to the fact that $\left( \begin{smallmatrix} \alpha_t & \beta_t \\ \gamma_t &
\delta_t \end{smallmatrix} \right)$
is $t$-independent modulo $M_2(B)$. $E \circ CH(\varphi) \in
\Ext^{-1/2}(T^2A,M_2(B) \otimes \mathbb K)$ is given by a continuous
approximate unit $\{u_{t}\}_{t \in [1,\infty)}$ in $B$ and a
sequence $t_0 < t_1 < t_2 < t_3 < \dots $ in $[1,\infty)$ such that
$$
E \circ CH(\varphi)(f \otimes z^k \otimes a) = q_{M_2(B \otimes
\mathbb K)}\left(  \left( \begin{matrix} \mathcal T^k & 0 \\ 0 &
\mathcal T^k \end{matrix}\right)  \sum_{n \in \mathbb Z}
\left( \begin{matrix}  f(u_n)\alpha_{t_n}(a) &  f(0)\beta_{t_n}(a) \\
f(0)\gamma_{t_n}(a)  & f(0)\delta_{t_n}(a)  \end{matrix} \right)
\otimes e_{n,n}\right) ,
$$
when $f \in T = \{ g \in C[0,1]: \ g(0) = g(1) \}$, $k \in
\mathbb Z$, $a \in A$. In this expression $t_n = |t_n|$ when
$n \leq 0$, $u_n = u_{t_n}$ when $n \geq 1$, and $u_n = 0$ when $n < 1$.

To describe $\Bott \circ E \circ CH(\varphi) \in \Ext^{-1/2}\left(A,M_2(B)
\otimes \mathbb K\right)$ we will ignore the $*$-isomorphism
$\hat{\Theta}$ appearing in the definition of $\Bott$. Then
$$
\Bott \circ E \circ CH(\varphi) \in \Ext^{-1/2}\left( A,M_4(B) \otimes
\mathbb K\right)
$$
is the difference $x -x_0$ of two elements, $x,x_0$, corresponding
to the projections $P$ and $P_0$, respectively. Using the explicit
description of $P$, we see that $x = [q_{M_4(B) \otimes \mathbb K}
\circ \psi]$, where $\psi : A \to M(M_4(B) \otimes \mathbb K)$ is
given, modulo $M_4(B) \otimes \mathbb K$, by the formula
\begin{eqnarray*}
\psi(a) &=&
\sum_{n \in \mathbb Z} \left( \begin{matrix}
s^2(u_n)\alpha_{t_n}(a) & 0 & sc_0(u_n)\alpha_{t_n}(a) & 0 \\
0 & 0 & 0 &  0  \\ sc_0(u_n)\alpha_{t_n}(a) & 0 &
[c_0^2 + c_1^2](u_n)\alpha_{t_n}(a) & \beta_{t_n}(a) \\
0 & 0 & \gamma_{t_n}(a) &  \delta_{t_n}(a) \end{matrix} \right)
\otimes e_{n,n} \\
&+&
\left(\begin{matrix} \mathcal T & 0 & 0 & 0 \\
0 & \mathcal T & 0 & 0 \\ 0 & 0 & \mathcal T^* & 0 \\
0 & 0 & 0 & \mathcal T^* \end{matrix} \right)
\sum_{n \in \mathbb Z} \left( \begin{matrix} 0 & 0 &
sc_1(u_n)\alpha_{t_n}(a) & 0 \\ 0 & 0 & 0 & 0  \\
sc_1(u_n)\alpha_{t_n}(a) & 0 & 0 & 0 \\ 0 & 0 & 0 & 0
\end{matrix} \right) \otimes e_{n,n} .
\end{eqnarray*}
$\psi$ is clearly unitarily equivalent to the map $\psi'$ given,
modulo $M_4(B) \otimes \mathbb K$, by
\begin{eqnarray*}
\psi'(a) &=&
\sum_{n \in \mathbb Z} \left( \begin{matrix}
s^2(u_n)\alpha_{t_n}(a) & sc_0(u_n)\alpha_{t_n}(a) & 0 & 0 \\
sc_0(u_n)\alpha_{t_n}(a)  &  [c_0^2 + c_1^2](u_n)\alpha_{t_n}(a) &
0 &  \beta_{t_n}(a)  \\0 & 0 & 0 & 0 \\ 0 & \gamma_{t_n}(a) & 0&
\delta_{t_n}(a) \end{matrix} \right) \otimes e_{n,n} \\
&+&
\left(\begin{matrix} \mathcal T & 0 & 0 & 0 \\
0 & \mathcal T^* & 0 & 0 \\ 0 & 0 & \mathcal T & 0 \\
0 & 0 & 0 & \mathcal T^* \end{matrix} \right)
\sum_{n \in \mathbb Z} \left( \begin{matrix} 0 &
sc_1(u_n)\alpha_{t_n}(a) &0 & 0 \\  sc_1(u_n)\alpha_{t_n}(a)
& 0 & 0 & 0  \\ 0 & 0 & 0 & 0 \\ 0 & 0 & 0 & 0 \end{matrix} \right)
\otimes e_{n,n} .
\end{eqnarray*}
Similarly, $x_0 = [q_{M_4(B) \otimes \mathbb K} \circ \psi_0]$,
where $\psi_0 : A \to M(M_4(B) \otimes \mathbb K)$ is given,
modulo $M_4(B) \otimes \mathbb K$, by the formula
\begin{equation*}
\psi_0(a) = \sum_{n \in \mathbb Z} \left( \begin{matrix} 0 & 0 & 0 & 0 \\
0 & 0 & 0 & 0  \\ 0 & 0 & \alpha_{t_n}(a) & \beta_{t_n}(a) \\
0 & 0 & \gamma_{t_n}(a) &  \delta_{t_n}(a) \end{matrix} \right)
\otimes e_{n,n}.
\end{equation*}
Define $\alpha,\beta,\gamma, \delta : A \to M(B \otimes \mathbb K)$ by
$$
\alpha(a) = \sum_{n \in \mathbb Z} \alpha_{t_n}(a) \otimes e_{n,n}, \
\beta(a) = \sum_{n \in \mathbb Z} \beta_{t_n}(a) \otimes e_{n,n}, \
\gamma(a) = \sum_{n \in \mathbb Z} \gamma_{t_n}(a) \otimes e_{n,n}, \
\delta(a) = \sum_{n \in \mathbb Z} \delta_{t_n}(a) \otimes e_{n,n},
$$
and $S,C_0,C_1 \in M(B \otimes \mathbb K)$ by
$$
S = \sum_{n \in \mathbb Z} s(u_n) \otimes e_{n,n}, \ C_0 =
\sum_{n \in \mathbb Z} c_0(u_n) \otimes e_{n,n}, \ C_1 =
\sum_{n \in \mathbb Z} c_1(u_n) \otimes e_{n,n}.
$$
Since $\alpha_{t_n}(a) - \alpha_{t_{n+1}}(a) \in B$ and
$\lim_{n \to \pm \infty} \alpha_{t_n}(a) - \alpha_{t_{n+1}}(a) = 0$
for all $a \in A$, we see that $[\alpha(a),\mathcal T] \in B \otimes
\mathbb K$ for all $a \in A$, i.e. $\mathcal T$ essentially commutes
with $\alpha$. The same is true, for the same reason, for $\beta$,
$\gamma$ and $\delta$. Similarly, we can assume that
$\lim_{n \to \infty} u_n - u_{n+1} = 0$, which implies
that also $S,C_0$ and $C_1$ essentially commute with $\mathcal T$.
Note that
$$
\psi'(a) = \left( \begin{matrix} S^2\alpha(a)  &
SC_1\mathcal T\alpha(a) + SC_0\alpha(a) & 0 & 0 \\
SC_1\mathcal T^*\alpha(a) + SC_0\alpha(a) &
\left(C_0^2 + C_1^2\right)\alpha(a) & 0 & \beta(a) \\
0 & 0 & 0 & 0 \\ 0 & \gamma(a) & 0 & \delta(a) \end{matrix} \right),
$$
while
$$
\psi_0(a) = \left( \begin{matrix} 0  & 0 & 0 & 0 \\
0 & 0 & 0 & 0 \\ 0 & 0 & \alpha(a) & \beta(a) \\ 0 & 0 & \gamma(a)  &
\delta(a) \end{matrix} \right).
$$
Set $\mathcal T_0 = \sum_{n \geq 1} e_{n,n-1} + \sum_{n \leq -1}
e_{n,n}$, and note that $\mathcal T_0 \in M(B \otimes \mathbb K)$
is an isometry such that $\mathcal T_0\mathcal T_0^* = 1 - e_{0,0}$.
Like $\mathcal T$, also $\mathcal T_0$ commutes with
$S,C_0,C_1,\alpha(a),\beta(a),\gamma(a)$ and $\delta(a)$,
modulo $B \otimes \mathbb K$. Set
$$
W_+ = \left( \begin{matrix} S & -C_0 - C_1\mathcal T \\
C_0  + C_1\mathcal T^* & S \end{matrix} \right) , \ W_- =
\left( \begin{matrix} e_{1,1} & - \mathcal T_0  \\
\mathcal T_0^*  & 0 \end{matrix} \right) \in
M_2\left(M(B \otimes \mathbb K)\right) .
$$
Then $W_-$ is a unitary while $W_+$ is unitary modulo
$M_2(B \otimes \mathbb K)$.

Furthermore, a calculation shows that
$$
\Ad \left( \begin{matrix} {W_+}^* & 0 \\ 0 & {W_-}^* \end{matrix} \right)
\circ \psi'(a) = \left( \begin{matrix} \alpha(a)  & 0 &
\left[C_0\mathcal T_0^* + C_1\mathcal T
\mathcal T_0^*\right]\beta(a) & 0 \\ 0 & 0 & S\mathcal T_0^*\beta(a) &
0 \\ \left[\mathcal T_0C_0  +
\mathcal T_0\mathcal T^*C_1\right]\gamma(a) & S\mathcal T_0\gamma(a) &
\mathcal T_0\mathcal T_0^*\delta(a) & 0 \\ 0 & 0 & 0 & 0 \end{matrix}
\right),
$$
modulo $M_4(B) \otimes \mathbb K$. Since $\lim_{n \to \infty}
\beta_{t_n}(a)s(u_n) = \lim_{n \to \infty} \gamma_{t_n}(a)s(u_n) = 0$,
we see that $S\beta(a),S\gamma(a) \in B \otimes \mathbb K$. It follows
that $ S\mathcal T_0^*\beta(a),  S\mathcal T_0\gamma(a) \in B \otimes
\mathbb K$. Similarly, since $\lim_{n \to \infty}
c_0(u_n)\beta_{t_n}(a) = \lim_{n \to \infty} c_0(u_n)\gamma_{t_n}(a) = 0$,
we find that $C_0\gamma(a) = P_-\gamma(a)$ and $C_0\beta(a) =
P_-\beta(a)$, modulo $B \otimes \mathbb K$, where $P_- =
\sum_{n \leq 0} e_{n,n}$. For a similar reason, we find that
$C_1\beta(a) = P_+\beta(a)$ and $C_1\gamma(a) = P_+\gamma(a)$,
modulo $B \otimes \mathbb K$, where $P_+  = \sum_{n \geq 1} e_{n,n}$.
It follows that $\left[C_0\mathcal T_0^* +
C_1\mathcal T\mathcal T_0^*\right]\beta(a) = P_-\beta(a) +
P_+\beta(a) = \beta(a)$ and $\left[\mathcal T_0C_0  +
\mathcal T_0\mathcal T^*C_1\right]\gamma(a) = P_-\gamma(a) +
P_+\gamma(a) = \gamma(a)$, modulo $B \otimes \mathbb K$.
Consequently,
$$
\Ad \left( \begin{matrix} {W_+}^* & 0 \\ 0 & {W_-}^* \end{matrix}
\right) \circ \psi'(a) = \left( \begin{matrix} \alpha(a)  & 0 &
\beta(a) & 0 \\ 0 & 0 & 0 & 0 \\ \gamma(a) & 0 &
\mathcal T_0\mathcal T_0^*\delta(a) & 0 \\ 0 & 0 & 0 & 0
\end{matrix} \right)
$$
modulo $M_4(B) \otimes \mathbb K$. Since conjugation by the unitary
$q_{M_2(B)}\left( \begin{smallmatrix} {W_+} &  \\  & {W_-} \end{smallmatrix}
\right)$ does not change the class in $\Ext^{-1/2}(A,M_4(B))$, we conclude that
$x = [q_{M_2(B)
\otimes \mathbb K} \circ \psi'']$, where
$$
\psi''(a) = \left( \begin{matrix} \alpha(a) &  \beta(a) \\
\gamma(a) & (1 -e_{0,0})\delta(a) \end{matrix} \right).
$$
Clearly, $x_0 = [q_{M_2(B) \otimes \mathbb K} \circ \psi_0'']$, where
$$
\psi_0''(a) = \left( \begin{matrix} \alpha(a) & \beta(a) \\
\gamma(a) & \delta(a) \end{matrix} \right).
$$
It follows that $x_0-x$ is represented by $a \mapsto q_{B \otimes
\mathbb K}\left(e_{0,0} \delta(a)\right)$, which represents
$e( - [\varphi])$ since $q_B \circ \delta$ represents $-[\varphi]$.
Hence $x - x_0 = e([\varphi])$. We have now shown that (\ref{inj})
commutes, and it follows that $CH$ is injective.

\subsection{$CH$ is surjective} Let $P,P_0 \in M_2(T^2)$ be the
projections used to define $\Bott$, cf. Section \ref{Bott-section}.
We can then define a unitary $U \in M_2(T^3) =
\left\{f \in C\left([0,1],M_2(C(\mathbb T^2))\right) :
f(0) = f(1) \right\}$ by
$$
U(s) = e^{2 \pi i s P}e^{-2 \pi i s P_0} =
\left(1 + (e^{2 \pi i s} -1)P\right)
\left(1 + (e^{-2 \pi i  s} -1)P_0\right).
$$
Note that $U(0) = U(1) = 1$. Since $P(0,z) = P(1,z) = P_0$ for all
$z \in \mathbb T$, it follows that
$$
U - 1 \in M_2(STS).
$$
Consequently the $*$-homomorphism $\phi_U : T \to M_2(T^3)$ given
by $\phi_U(f) = f(U)$ has the property that
$\phi_U(S) \subseteq M_2(STS)$. We can therefore define a
map $\bold B : [[T^3A,STSA;B]] \to [[TA,SA;B]]$ such that
$\bold B[\varphi] = [\varphi']$, where
$$
\varphi'_t(a) = \tilde{\Theta}\circ \left( \id_{M_2} \otimes
\varphi_t\right)\left((\phi_U \otimes \id_A)(a)\right)
$$
$a \in TA$. To compare $\bold B$ with the other maps we have
in play, let $j : [[T^3A,ST^2A;B \otimes \mathbb K]] \to
[[T^3A,STSA;B \otimes \mathbb K]]$ be the forgetfull
homomorphism obtained from the fact that $STSA \subseteq ST^2A$.
We claim that the diagram
\begin{equation}\label{bold}
\begin{xymatrix}{
[[T^3A,ST^2A;B]] \ar[d]_-j \ar[r]^-{\Bott} & [[TA,SA;B]]\\
[[T^3A,STSA;B]] \ar[ur]_-{\bold B} & {} }
\end{xymatrix}
\end{equation}
commutes.  To see this define $\pi_1, \lambda_1 : T \to M_2(T^3)$ by
$\pi_1(f) = f P$ and $\lambda_1(f) = fP + f(0)(1-P)$. Similarly, we
set $\pi_0(f) = fP_0$ and $\lambda_0(f) = fP_0 + f(0)(1-P_0)$.
Note that all these $*$-homomorphisms take $S$ into $M_2(ST^2)$.
Let $z \in T$ denote the canonical unitary generator; the identity
function on $\mathbb T$. For $\theta \in [0,\frac{\pi}{2}]$, set
$$
V_{\theta} = \left( \begin{smallmatrix} zP + 1 - P & {} \\ {} &
1 \end{smallmatrix} \right)
\left( \begin{smallmatrix} \cos \theta  & \sin \theta  \\ -\sin
\theta & \cos \theta \end{smallmatrix} \right)
\left( \begin{smallmatrix} z^*P_0 + 1 - P_0 & {} \\ {} &
1 \end{smallmatrix} \right)\left( \begin{smallmatrix}
\cos \theta  & -\sin \theta  \\ \sin \theta & \cos \theta
\end{smallmatrix} \right)\left( \begin{smallmatrix} 1  & {} \\
{} & zP_0 + 1 - P_0 \end{smallmatrix} \right),
$$
which gives a homotopy of unitaries in $M_4(T^3)$ connecting
$\left( \begin{smallmatrix} U  & {} \\ {} & zP_0 + 1 - P_0
\end{smallmatrix} \right)$ to $ \left( \begin{smallmatrix}
zP + 1 - P & {} \\ {} & 1 \end{smallmatrix} \right)$. Note
that when we substitute $1$ for $z$ in the formula for
$V_{\theta}$, we get $1$ for each $\theta$. It follows that
$\phi_U \oplus \lambda_0$ is homotopic to $\lambda_1 \oplus \ev$,
where $\ev(f) = f(0)$, via a path of $*$-homomorphisms taking
$S$ into $M_2(ST^2)$. Since $[\left(\id_{M_2} \otimes
\varphi\right) \circ \lambda_i] = [\left(\id_{M_2} \otimes
\varphi\right) \circ \pi_i], i = 0,1$, in $[[TA,SA;B]]$
for all $\varphi$ by Lemma \ref{triv}, we conclude that
$$
\tilde{\Theta}\circ \left( \id_{M_2} \otimes \varphi_t\right)
\left((\phi_U \otimes \id_A)(a)\right) \oplus \tilde{\Theta} \circ
\left(\id_{M_2} \otimes \varphi_t\right)\left(P_0 \otimes a\right)
$$
defines the same element of $[[TA,SA;B]]$ as $\tilde{\Theta}
\left((\id_{M_2} \otimes \varphi_t)\left(P \otimes a\right)\right)$.
This establishes the commutativity of (\ref{bold}).

Let $\varphi = \left( \varphi_t\right)_{t \in [1,\infty)} :
TA \to M(B)$ be an asymptotic homomorphism, constantly extended
from $SA$. Let $\left(\varphi_{t_n}\right)_{n \in \mathbb N}$
be a discretization of $\varphi$. For each $a \in TA, t \in
[1,\infty)$, set
$$
\overline{\varphi}_t(a) = \sum_{n \geq 1} \varphi_{\max\{t_n,t\}}(a)
\otimes e_{n,n}  +\sum_{n \leq 0} \varphi_{\max\{t_{|n|},t\}}(a)
\otimes e_{n,n},
$$
which is an element of $M(B \otimes \mathbb K)$. Then
$\overline{\varphi} =
\left( \overline{\varphi}_t\right)_{t \in [1,\infty)}$
is an asymptotic homomorphism which is essentially constant,
i.e. $\overline{\varphi}_t(a) - \overline{\varphi}_s(a) \in
B \otimes \mathbb K$ for all $t,s \in [1,\infty)$.
Furthermore, $\overline{\varphi}_t(a)$ commutes with the
two-sided shift $\mathcal T$ modulo $B \otimes \mathbb K$.
For each $n \geq 1$, set
$$
v_n = \sum_{i \in \mathbb Z} v_n^{(i)}\otimes e_{i,i},
$$
where
$$
v_n^{(i)} = \begin{cases} 1 , & i \leq 0 \\ \frac{n-i}{n}, &
\ 1 \leq i \leq n, \\ 0, & i \geq n  . \end{cases}
$$
Then set
\begin{equation}\label{ogg}
v_t = (t - n)v_{n+1} + (n+1 -t)v_n,
\end{equation}
$t  \in [n,n+1]$. It follows that $[v_t,\overline{\varphi}_s(a)] = 0$
for all $a,s,t$, and that $\lim_{t \to \infty} [\mathcal T,f(v_t)] = 0$
for all $f \in C[0,1]$ for which $f(0) = f(1)$. We can therefore
define an asymptotic homomorphism $\beta(\varphi) : T^3A \to
M(B \otimes \mathbb K)$ determined, up to asymptotic equality,
by the condition that
$$
\lim_{t \to \infty} \beta(\varphi)_t(f \otimes z^k \otimes a) -
\overline{\varphi}_t(a)f(v_t)\mathcal T^k = 0
$$
when $f \in C(\mathbb T), k \in \mathbb Z$ and $a \in TA$. Since
$\overline{\varphi}_t(a)f(v_t)\mathcal T^k -
\overline{\varphi}_s(a)f(v_s)\mathcal T^k \in B \otimes \mathbb K$
for all $s,t,a,f,k$, and $\overline{\varphi}_t(a)f(v_t)
\mathcal T^k \in B \otimes \mathbb K$, when $a \in SA$ and
$f \in S$, we can arrange that $\beta(\varphi)$ is essentially
constant and that
$$
\beta(\varphi)_t(STSA) \subseteq B \otimes \mathbb K,
$$
for all $t \in [1,\infty)$, cf. the construction in Remark
\ref{remark}. We get in this way a map $\beta : [[TA,SA;B]]
\to [[T^3A,STSA;B \otimes \mathbb K]]$ such that
$\beta[\varphi] = [\beta(\varphi)]$. We claim that the diagram
\begin{equation}\label{natdiag2}
\begin{xymatrix}{
[[TA,SA;B]] \ar[r]^-E \ar@/_1pc/[ddr]_-{\beta} &
\Ext^{-1/2}(T^2A, B \otimes \mathbb K) \ar[d]^-{CH} \ar[r]^-{\Bott} &
\Ext^{-1/2}(A,B\otimes \mathbb K) \ar[d]^-{CH} \\
& [[T^3A,ST^2A;B\otimes \mathbb K]] \ar[r]^-{\Bott}
\ar[d]^-j & [[TA,SA; B\otimes \mathbb K]] \\
& [[T^3A,STSA;B \otimes \mathbb K]] \ar[ur]_-{\bold B} }
\end{xymatrix}
\end{equation}
commutes. Since the square commutes by the naturality of the
extended Connes-Higson construction, cf. (\ref{natdiag}), and the right
triangle commutes by (\ref{bold}), it suffices to show that the left
triangle commutes, i.e. we must show that $\beta = j \circ CH \circ E$.
Let therefore $\varphi : TA \to M(B)$ be a constantly extended
asymptotic homomorphism. $E(\varphi)$ is given by (\ref{klaus1})
and (\ref{klaus2}) for an appropriate discretization
$\left( \varphi_{t_n} \right)_{n \in \mathbb N}$ of $\varphi$,
and the inverse $-E(\varphi)$ is given by (\ref{klaus3}) and
(\ref{klaus4}). We shall use the constructions of
Remark \ref{remark} in order to get a workable description of
$CH \circ E(\varphi)$. In particular, we refer to Remark \ref{remark}
for the notation used in the following. Let $\{u_t\}_{t \in [1,\infty)}$
be a continuous approximate unit in $B \otimes \mathbb K$ satisfying
the requirements needed to define $CH \circ E(\varphi)$, cf.
(\ref{CH1})-(\ref{CH3}). Then $CH \circ E[\varphi]$ is represented
by an asymptotic homomorphism $\tilde{\Theta} \circ \psi$, where the
asymptotic homomorphism $\psi : T^3A \to M_2(M(B\otimes \mathbb K))$
satisfies that $\psi_t(h \otimes z^k \otimes f)$ asymptotically
agrees with
\begin{equation}\label{endelig1}
 \left( \begin{smallmatrix} h(u_t) & {} \\ {}  & h(1) \end{smallmatrix}
\right )\left( \begin{smallmatrix} \mathcal T^k & {} \\ {}  &
\mathcal T^k \end{smallmatrix} \right ) U_{t} \Gamma_t(f)U_{t}^* ,
 \end{equation}
for all $h \in T, k \in \mathbb Z$ and $f \in TA$. Note that
\begin{equation*}
\begin{split}
&U_{t} \Gamma_t(f) U_{t}^* = \\
& \sum_{i \in \mathbb Z} \left( \begin{matrix} \alpha_i(t)
\varphi^i_t(f) + (1-\alpha_i(t))\varphi^i_t(f(0)) &
\sqrt{\alpha_i(t) - \alpha_i(t)^2}\left( \varphi^i_t(f) -
\varphi^i_t(f(0))\right) \\ \sqrt{\alpha_i(t) - \alpha_i(t)^2}
\left( \varphi^i_t(f) - \varphi^i_t(f(0))\right) &
(1 -\alpha_i(t))\varphi^i_t(f) + \alpha_i(t)\varphi^i_t(f(0))
\end{matrix} \right) \otimes e_{i,i} .
\end{split}
\end{equation*}
Let $F_1 \subseteq F_2 \subseteq F_3 \subseteq \dots$ be a sequence
of finite sets with dense union in $TA$.
Let $v_n, n \in \mathbb N$, be an approximate unit in $B$ such that
$$
\left\|\left[v_n, \varphi_t(f)\right]\right\| \leq \frac{1}{n}
$$
and
$$
\left\| \left( v_n - 1 \right) \left( \varphi_t(f) -
\varphi_t(f(0))\right)\right\| \leq \frac{1}{n}
$$
for all $t \in [1,3n], f \in F_n$. Set
$$
\beta_i = \begin{cases} 1 , & \ i \in \{-n,i-n+1, \dots, n-1,n\} \\
\frac{2n -i}{n}, & i \in \{n+1,n+2, \dots, 2n\} \\  \frac{2n+i}{n}, &
\ i \in \{-2n-1, -2n+1, \dots , -n-1\} \\ 0, & |i| > 2n, \end{cases}
$$
$\tilde{v}_n = \sum_{i \in \mathbb Z} \beta_i v_n \otimes e_{i,i}$, and
$$
w_t = (n+1 -t)\tilde{v}_{n} + (n -t)\tilde{v}_{n+1},
$$
when $t \in [n,n+1]$. Then $\{w_t\}_{t \in [0,\infty)}$ is a
continuous approximate unit in $B \otimes \mathbb K$ such that
the requirements needed to define $CH \circ E(\varphi)$, cf.
(\ref{CH1})-(\ref{CH3})
hold for $w_t$ in place of $u_t$. We can therefore
work with this path instead of $\{u_t\}_{t \in [1,\infty)}$ in
(\ref{endelig1}).
Thanks to (\ref{ogg8}), (\ref{endelig1}) then becomes asymptotically the same as
\begin{equation}\label{endelig11}
  U_{t}  \left[U_{t}^* \left( \begin{smallmatrix} h(w_t) & {} \\
{}  & h(1) \end{smallmatrix} \right ) U_{t}\right]
\left( \begin{smallmatrix} \mathcal T^k & {} \\ {}  &
\mathcal T^k \end{smallmatrix} \right )\left( \sum_{i \in
\mathbb Z}\left( \begin{smallmatrix} \varphi^i_t(f) & {} \\
{}  & \varphi^i_t(f(0)) \end{smallmatrix} \right ) \otimes
e_{i,i} \right)U_{t}^* .
\end{equation}
Note that conjugation by $U_{t}$ induces the identity map in
$[[T^3A,ST^2A; M_2(B \otimes \mathbb K)]]$. We see therefore
from (\ref{endelig11}) that $CH \circ E[\varphi]$ is represented
by an asymptotic homomorphism $\tilde{\Theta} \circ \psi'$,
such that $\psi'_t(h \otimes z^k \otimes f)$ asymptotically agrees with
$$
 h(Y_t) \left( \sum_{i \in \mathbb Z} \left( \begin{smallmatrix} 1 & {} \\
{} & 1 \end{smallmatrix} \right) \otimes e_{i,i-1}\right)^k
\left(\sum_{i \in \mathbb Z}\left( \begin{smallmatrix}
\varphi^i_t(f) & {} \\ {}  & \varphi^i_t(f(0)) \end{smallmatrix}
\right )\otimes e_{i,i}\right),
$$
where
$$
Y_t = U_{t}^* \left( \begin{smallmatrix} w_t & {} \\ {}  &
1 \end{smallmatrix} \right ) U_{t}.
$$
We want to substitute $Y_t$ with something else. To this end write
$$
\overline{\varphi}_t(f) = \sum_{i \in \mathbb Z} \varphi^i_t(f)
\otimes e_{i,i} ,
$$
and let
$$
Y_t = \left( \begin{smallmatrix} Y^{11}_t & Y^{12}_t \\
Y^{21}_t & Y^{22}_t \end{smallmatrix} \right)
$$
be the $2 \times 2$-matrix decomposition of $Y_t$. Set $Q =
\sum_{i \leq 0} e_{i,i}$. The significant properties of $Y_t$ are
the following:
\begin{enumerate}
\item[0)] $0 \leq Y_t \leq 1$,
\item[1)] $\lim_{t \to \infty} [Y^{11}_t,\overline{\varphi}_t(f)] =
\lim_{t \to \infty} [Y^{22}_t,\overline{\varphi}_t(f(0))] = 0$ for
all $f \in TA$,
\item[2)] $\lim_{t \to \infty} Y^{12}_t
\left( \overline{\varphi}_t(f) - \overline{\varphi}_t(f(0))\right) =
\lim_{t \to \infty} Y^{21}_t \left( \overline{\varphi}_t(f) -
\overline{\varphi}_t(f(0))\right) = 0$  for all $f \in TA$,
\item[3)] $\lim_{t \to \infty} \left[ Y_t,
\left( \begin{smallmatrix} \mathcal T & {} \\ {}  & \mathcal T
\end{smallmatrix} \right) \right] = 0$,
\item[4)] $Y_t \left( \begin{smallmatrix} \overline{\varphi}_t(f) & {} \\
{}  & 0 \end{smallmatrix} \right) =  \left( \begin{smallmatrix}
Q \overline{\varphi}_t(f) & {} \\ {}  &  0 \end{smallmatrix}
\right)$  modulo $M_2(B \otimes \mathbb K)$ for all $t$ and all
$f \in SA$,
\end{enumerate}
which are all easy to check. Note that 4) implies
\begin{enumerate}
\item[5)] $\left[ Y_t,  \left( \begin{smallmatrix} \mathcal T &
{} \\ {}  & \mathcal T \end{smallmatrix} \right) \right]
\left( \begin{smallmatrix}  \overline{\varphi}_t(f) & {} \\
{}  &  0 \end{smallmatrix} \right) \in M_2(B \otimes \mathbb K)$
for all $t$ and all $f \in SA$,
\end{enumerate}
since
$$
\left[  \left( \begin{smallmatrix} \mathcal T & {} \\
{}  & \mathcal T \end{smallmatrix} \right),\left( \begin{smallmatrix}
\overline{\varphi}_t(f) & {} \\ {}  &  0 \end{smallmatrix} \right)
\right], \ \left[  \left( \begin{smallmatrix} \mathcal T & {} \\
{}  & \mathcal T \end{smallmatrix} \right),\left( \begin{smallmatrix}
Q \overline{\varphi}_t(f) & {} \\ {}  &  0 \end{smallmatrix} \right)
\right] \in M_2(B \otimes \mathbb K)
$$
for all $t$ and all $f \in SA$.

Put
$$
Y^{\lambda}_t = (1-\lambda) Y_t + \lambda \left( \begin{smallmatrix}
v_t & {} \\ {} & 0 \end{smallmatrix} \right),
$$
where $v_t$ is defined by (\ref{ogg}). Then $Y^{\lambda}_t$
satisfies 0)-4) for all $\lambda \in [0,1]$. It follows from 0)-3)
that we can define an asymptotic homomorphism $\Phi : T^3A \to
M_2(M(IB \otimes \mathbb K))$ such that $\Phi_t(h \otimes z^k
\otimes f)(\lambda), \lambda \in [0,1]$, asymptotically agrees with
$$
h(Y^{\lambda}_t) \left( \begin{smallmatrix} \mathcal T^k & {} \\
{}  & \mathcal T^k \end{smallmatrix} \right)
\left( \begin{smallmatrix} \overline{\varphi}_t(f) & {} \\
{}  & \overline{\varphi}_t(f(0)) \end{smallmatrix} \right)
$$
for all $h \in T, k \in \mathbb Z$ and $f \in TA$. We claim
that we can arrange that
$$
\Phi_t(STSA) \subseteq M_2(IB \otimes \mathbb K)
$$
for all $t$. Since $S = \{\mu \in C[0,1] : \mu(1) = \mu(0) = 0\}$
is generated by the function $s \mapsto e^{2 \pi is} -1$, it
suffices for this purpose to check that
$$
   \left( e^{2 \pi i  \pm Y^{\lambda}_t} - 1\right)
\left( \begin{smallmatrix} \mathcal T^k & {} \\ {}  &
\mathcal T^k \end{smallmatrix} \right) \left( \begin{smallmatrix}
\overline{\varphi}_t(f) & {} \\ {}  & \overline{\varphi}_t(f(0))
\end{smallmatrix} \right)$$
is in $ M_2(B \otimes \mathbb K)$ for $\lambda \in [0,1]$,
$k \in \mathbb Z$ and $f \in SA$.
This follows from 5) and 4) because we see that
\begin{equation*}
\begin{split}
&\left( e^{2 \pi i  \pm Y^{\lambda}_t} - 1\right)
\left( \begin{smallmatrix} \mathcal T^k & {} \\ {}  &
\mathcal T^k \end{smallmatrix} \right) \left( \begin{smallmatrix}
\overline{\varphi}_t(f) & {} \\ {}  & \overline{\varphi}_t(f(0))
\end{smallmatrix} \right) \\
& =  \left( \begin{smallmatrix} \mathcal T^k & {} \\ {}  &
\mathcal T^k \end{smallmatrix} \right) \left( e^{2 \pi i
\pm Y^{\lambda}_t} - 1\right) \left( \begin{smallmatrix}
\overline{\varphi}_t(f) & {} \\ {}  & \overline{\varphi}_t(f(0))
\end{smallmatrix} \right) \\
&  = \left( \begin{smallmatrix} \mathcal T^k & {} \\ {}  &
\mathcal T^k \end{smallmatrix} \right) \left( \begin{smallmatrix}
\left(e^{2 \pi i \pm Q} -1\right)\overline{\varphi}_t(f) & {} \\
{}  & 0 \end{smallmatrix} \right)  = 0,
\end{split}
\end{equation*}
modulo $M_2(B \otimes \mathbb K)$. Thus $\Phi$ gives us a homotopy
of, not necessarily constantly, extended asymptotic homomorphisms. At both
ends the asymptotic homomorphisms are constantly extended so we
can conclude from Theorem \ref{Thm3} that $j \circ CH \circ E[\varphi]
\in [[T^3A,STSA;B\otimes \mathbb K]]$ is represented by an
asymptotic homomorphism $\psi$ such that $\psi_t(h \otimes z^k
\otimes f)$ essentially (i.e. modulo $B\otimes\mathbb K$)
and asymptotically agrees with
$$
h(v_t) \mathcal T^k \overline{\varphi}_t(f)
$$
when $h \in T, k \in \mathbb Z$ and $f \in TA$. But this is
$\beta(\varphi)$, so we have shown that the diagram
(\ref{natdiag2}) commutes.

It suffices now to show that $\bold B \circ \beta = e$. To this
end define for each $f \in TA$ an element $H(f) \in M(IB \otimes
\mathbb K)$ such that
$$
H(f)(\lambda) = \sum_{i \in \mathbb Z} \varphi_{\lambda
\max\{t_{|i|},t\} + (1-\lambda)t}(f)\otimes e_{i,i} .
$$
We can then define an asymptotic homomorphism $\Psi :
T^3A \to M(IB \otimes \mathbb K)$ such that $\Psi_t(h
\otimes z^k \otimes f)$ asymptotically agrees with
$$
h(v_t)\mathcal T^k H(f) .
$$
Since $h(v_t)\mathcal T^k H(f) \in IB \otimes \mathbb K$,
when $h \in S,f \in SA$, we get a homotopy of (typically not
constantly) extended asymptotic homomorphisms showing that
$\beta(\varphi)$ defines the same element in $[[T^3A,STSA;
B \otimes \mathbb K]]$ as an asymptotic homomorphism $\psi$
with the property that $\psi_t(h \otimes z^k \otimes f)$
asymptotically agrees with
$$
h(v_t)\mathcal T^k \sum_{i \in \mathbb Z} \varphi_t(f) \otimes e_{i,i}
$$
for all $h \in T, k \in \mathbb Z, f \in TA$. To compare this
with $\varphi$, define an asymptotic homomorphism $\varphi
\otimes \id_{\mathbb K^+} : TA \otimes \mathbb K^+ \to M(B
\otimes \mathbb K)$ such that $\left(\varphi \otimes
\id_{\mathbb K^+}\right)_t(f \otimes x)$ asymptotically agrees with
$$
\iota\left(\varphi_t(f) \otimes x\right),
$$
where $\iota : M(B) \otimes \mathbb K^+ \to M(B \otimes \mathbb K)$
is the canonical embedding. Since $\varphi_t(SA) \subseteq B$ we
can arrange that $\left( \varphi \otimes
\id_{\mathbb K^+}\right)_t\left(SA \otimes \mathbb K\right)
\subseteq B \otimes \mathbb K$ for all $t$. Define also a
continuous path $A_t, t \in [1,\infty)$, of contractions in
$M_2(T\mathbb K^+)$ by
$$
A_t = \left[zQ_t + 1 - Q_t\right]\left[z^*P_0 - 1 - P_0\right],
$$
where $z \in T$ is the identity function, $P_0 =
\left( \begin{smallmatrix} 0 & {} \\ {} & 1 \end{smallmatrix}
\right)$ and
$$
Q_t = \left( \begin{matrix} s^2(v_t) & sc_0(v_t) + sc_1(v_t)\mathcal T \\
\mathcal T^* sc_1(v_t) + sc_0(v_t) & 1 - s^2(v_t) \end{matrix} \right) .
$$
Then, by definition, $\bold B(\psi) = \tilde{\Theta} \circ \psi'$, where
$\psi' : TA \to M_2(M(B \otimes \mathbb K))$ is an asymptotic
homomorphism such that $\psi'_t(z^k \otimes a)$ asymptotically
agrees with
$$
\left(\id_{M_2} \otimes \left(\varphi \otimes \id_{\mathbb K^+}
\right)\right)_t \left( A_t^k \otimes a\right)
$$
for all $k \in \mathbb Z, a \in A$. Since $v_t$ asymptotically
commutes with $\mathcal T$, we see that $\lim_{t \to \infty}
\left\|Q_t^2 - Q_t\right\| = 0$. Hence a standard application
of spectral theory gives us a continuous path $\{P_t\}_{t \in
[1,\infty)}$ of projections in $M_2\left(\mathbb K^+\right)$ such that
$$
\lim_{t \to \infty} \left\|P_t - Q_t\right\| = 0.
$$
Since $Q_t - P_0 \in M_2\left(\mathbb K\right)$ we can arrange (or rather, the
standard procedure will automatically ensure) that
$$
P_t - P_0 \in M_2\left(\mathbb K\right) .
$$
It follows that
$$
U_t =   \left[zP_t + 1 - P_t\right]\left[z^*P_0 - 1 - P_0\right]
$$
is a continuous path of unitaries in $M_2(T\mathbb K^+)$ such that
\begin{equation}\label{fredag1}
\lim_{t \to \infty} \left\|U_t - A_t \right\| = 0
\end{equation}
and
$$
U_t - 1 \in M_2(S\mathbb K).
$$
Note that $\psi'$ is then asymptotically equivalent
to the asymptotic homomorphism
$ \left(\id_{M_2} \otimes \left(\varphi \otimes \id_{\mathbb K^+}
\right)\right) \circ \phi_{U_t}$, where $\phi_{U_t} :  TA \to M_2\left(TA
\otimes \mathbb K^+\right)$ is the family of
$*$-homomorphisms defined in such a way that
$\phi_{U_t}(z^k \otimes a) = U_t^k \otimes a$.
If $V \in M_2\left(T\mathbb K^+\right)$ is any other
unitary which is homotopic to $U_t$ for any sufficiently large $t$
within the subgroup of the
unitary group of $M_2\left(T\mathbb K^+\right)$ consisting of the unitaries $W$
such that $W - 1 \in M_2(S\mathbb K)$, then $\psi'$ is homotopic
to  $ \left(\id_{M_2} \otimes \left(\varphi \otimes
\id_{\mathbb K^+}\right)\right) \circ \phi_{V}$. Now note that
by definition $U_t$ is the image of the projction $P_t$ under
the loop-construction implementing the Bott-isomorphism
$K_0(\mathbb K) \to K_1(S\mathbb K)$. It is easy to see
that $P_t$, for all large $t$, represents the generator
$1$ under the canonical isomorphism $K_0(\mathbb K) \simeq
\mathbb Z$, and it follows from this that $U_t$ is homotopic,
within the indicated subgroup of the unitary group of $M_2\left(T \mathbb K^+\right)$,
to the unitary
$$
R = \left( \begin{matrix} ze_{00} + \sum_{i \in \mathbb Z \backslash \{ 0\}} e_{ii} & 0 \\ 0 &  \sum_{i \in \mathbb Z} e_{ii} \end{matrix} \right) .
$$
Consequently $[\psi'] = e[\varphi] + [\varphi_0]$ in $[[TA,SA;B \otimes \mathbb K]]$, where $\varphi_0:
TA\to M(B \otimes \mathbb K)$ is an asymptotic homomorphism which factors
through the evaluation map $\ev : TA \to A$. $\varphi_0$ represents
zero in $[[TA,SA;B ]]$ by Lemma \ref{triv} and we conclude that $[\psi'] = e[\varphi]$.

\section{Conclusion} Our main result, Theorem \ref{mainth}, shows that the map $\Ext^{-1/2}(A,B) \to [[SA,B]]$ arising from the Connes-Higson construction as defined in \cite{CH} factors as
\begin{equation*}
\begin{xymatrix}{
\Ext^{-1/2}(A,B) \ar[r] \ar[d]  & [[TA,SA;B]] \ar[dl] \\
[[SA,B]] &  }
\end{xymatrix}
\end{equation*}
such that the horizontal map is an isomorphism. Thus the question whether or not the vertical $CH$-map is an isomorphism has been transformed to a question which solely involves homotopy classes of asymptotic homomorphisms. Specifically the question is now whether or not the restriction map $  [[TA,SA;B]] \to [[SA,B]]$ is an isomorphism. It is with some regret that we must report that we haven't been able to decide the latter.

\end{document}